\documentclass{gtmon_a}
\pdfoutput=1

\usepackage{pinlabel}
\input{diagxy}
\let\to\rightarrow
\let\barrsquare\square
\let\square\undefined

\proceedingstitle{Heegaard splittings of 3--manifolds (Haifa 2005)}
\conferencestart{10 July 2005}
\conferenceend{19 July 2005}
\conferencename{Heegaard splittings of 3--manifolds}
\conferencelocation{Haifa}

\editor{Cameron Gordon}
\givenname{Cameron}
\surname{Gordon}

\editor{Yoav}
\givenname{Yoav}
\surname{Moriah}

\title{The Heegaard structure of Dehn filled manifolds}

\author{Yoav Moriah}
\givenname{Yoav}
\surname{Moriah}
\address{Department of Mathematics\\
Technion -- Israel Institute of Technology\\\newline
Haifa 32000\\
Israel}
\email{ymoriah@tx.technion.ac.il}
\urladdr{}

\author{Eric Sedgwick}
\givenname{Eric}
\surname{Sedgwick}
\address{DePaul CTI\\
243 S Wabash Avenue\\
Chicago IL 60604\\
USA}
\email{esedgwick@cs.depaul.edu}
\urladdr{}

\volumenumber{12}
\issuenumber{}
\publicationyear{2007}
\papernumber{08}
\startpage{233}
\endpage{263}

\doi{}
\MR{}
\Zbl{}

\arxivreference{0706.1927}


\received{23 May 2007}
\revised{3 July 2007}
\accepted{3 July 2007}
\published{3 December 2007}
\publishedonline{3 December 2007}
\proposed{}
\seconded{}
\corresponding{}
\version{}

\makeatletter
\def\cnewtheorem#1[#2]#3{\newtheorem{#1}{#3}[section]
\expandafter\let\csname c@#1\endcsname\c@pro}

\let\xysavmatrix\xymatrix
\def\xymatrix{\disablesubscriptcorrection\xysavmatrix}
\AtBeginDocument{\let\bar\wbar\let\tilde\wtilde\let\hat\what}

\makeop{lcm}

\cnewtheorem{thm}[pro]{Theorem}
\cnewtheorem{lem}[pro]{Lemma}
\cnewtheorem{clm}[pro]{Claim}
\cnewtheorem{example}[pro]{Example}
\cnewtheorem{cnj}[pro]{Conjecture}
\cnewtheorem{cor}[pro]{Corollary}
\cnewtheorem{quest}[pro]{Question}
\cnewtheorem{prob}[pro]{Problem}

\theoremstyle{definition}
\cnewtheorem{dfn}[pro]{Definition}
\cnewtheorem{exer}[pro]{Exercise}

\makeatother  
\makeautorefname{pro}{Proposition}

\theoremstyle{remark}
\newtheorem*{rmk}{Remark}

\begin{document}

\begin{htmlabstract}
We expect manifolds obtained by Dehn filling to inherit properties from
the knot manifold.  To what extent does that hold true for the Heegaard
structure? We  study four changes to the Heegaard structure that may occur
after filling: (1) Heegaard genus decreases, (2) a new Heegaard surface
is created, (3) a non-stabilized Heegaard surface destabilizes, and (4)
two or more non-isotopic Heegaard surfaces become isotopic.
We survey general results that give quite satisfactory restrictions
to phenomena (1) and (2) and, in a parallel thread, give a complete
classification of when all four phenomena occur when filling most
torus knot exteriors. This latter thread yields sufficient (and perhaps
necessary) conditions for the occurrence of phenomena (3) and (4).
\end{htmlabstract}

\begin{abstract}
We expect manifolds obtained by Dehn filling to inherit properties from
the knot manifold.  To what extent does that hold true for the Heegaard
structure? We  study four changes to the Heegaard structure that may occur
after filling: (1) Heegaard genus decreases, (2) a new Heegaard surface
is created, (3) a non-stabilized Heegaard surface destabilizes, and (4)
two or more non-isotopic Heegaard surfaces become isotopic.
We survey general results that give quite satisfactory restrictions
to phenomena (1) and (2) and, in a parallel thread, give a complete
classification of when all four phenomena occur when filling most
torus knot exteriors. This latter thread yields sufficient (and perhaps
necessary) conditions for the occurrence of phenomena (3) and (4).
\end{abstract}
\maketitle

\section{Introduction}

Let $X$ be a {\it knot manifold}, that is a compact, orientable and
irreducible 3--manifold with a single torus boundary component.
There are many results demonstrating that most of the manifolds
obtained by filling inherit properties from the knot manifold. We
would also expect the Heegaard structure of filled manifolds to be
closely related to the Heegaard structure of the knot manifold. For
example, it is easy to see that every Heegaard surface for the knot
manifold is a Heegaard surface for each filled manifold.  In
particular, this implies that the Heegaard genus of $X$ is an upper
bound on the genus of each filled manifold. However, the Heegaard
structure of a filled manifold can differ from that of the knot
manifold. Here are four ways that this could occur:

\begin{enumerate}
\item Heegaard genus decreases.
\item A new Heegaard surface is created.
\item A non-stabilized Heegaard surface destabilizes.
\item Two or more non-isotopic Heegaard surfaces become isotopic.
\end{enumerate}

 By a {\it new} Heegaard surface, we mean that a filled manifold
contains a Heegaard surface that is not isotopic (in the filled
manifold) to a Heegaard surface for the knot manifold $X$. When the
genus decreases (1), the filled manifold has a Heegaard surface of
lower genus than every Heegaard surface for $X$.  Indeed, it is a
new Heegaard surface. So, restricting (2) also restricts (1).

In each of these cases, we would like to either demonstrate that the
set of fillings for which the phenomenon occurs is special, for example
finite, a line of slopes, and/or conclude that the Heegaard
surface(s) in question are special in some regard, for example
$\gamma$--primitive, padded, or boundary stabilized.

In \fullref{secResults} we survey known work that gives quite
satisfactory restrictions to phenomena (1) and (2). We also give an
extended example: Dehn filling on a torus knot exterior, for which
we have almost complete knowledge.  We are able to completely
specify the fillings for which each of these four phenomena occur.
This also illustrates sufficient conditions for (3) and (4) to
occur.

\subsection*{Acknowledgements}
We would like to thank Ryan Derby-Talbot, David Bachman, Jesse Johnson
and Jennifer Schultens for helpful discussions. Also thanks to DePaul
CTI and the Mathematics Department of the Technion for their hospitality.

This research was supported by Grant No.~2002039 from the US--Israel
Binational Science Foundation (BSF), Jerusalem, Israel.

\section{Background}

\subsection[]{Dehn filling and slopes}

For simplicity assume that $X$ is a {\it knot manifold}, an
orientable, irreducible 3--manifold with boundary consisting of a
single incompressible torus.   Much of the discussion, and some of
the results quoted, also pertain to manifolds with multiple torus
boundary components but it will simplify our discussion not to
consider them here.

A {\it slope} $\alpha$ is the isotopy class of a simple closed curve
$\alpha$ in the boundary torus $\partial X$.   With a choice of
basis, for example the meridian longitude pair $(\mu,\lambda)$ for a
knot exterior in $S^3$, we can naturally identify the set of slopes
with $\mathbb{Q} \cup \{\infty = 1/0\}$. It will be important to be
able to identify $L_\alpha$, the {\it line of slopes associated with
a given slope $\alpha$}. These are precisely the slopes that meet
$\alpha$ exactly once: $L_\alpha = \{ \beta ~|~
\Delta(\alpha,\beta)=1 \}$, where $\Delta(\alpha,\beta)$ indicates
the geometric intersection number.   We will also construct a {\it
line of lines associated with $\alpha$}, $LL_\alpha$, which is the
set of slopes $\gamma \in L_\beta$ for some $\beta \in L_\alpha$,
that is, $LL_\alpha = \{ \gamma ~|~ \exists \beta~ s.t.
~\Delta(\alpha,\beta) = \Delta(\beta,\gamma)=1 \}$.

\subsection{Triviality conditions}
\label{secTrivConditions}

We assume that the reader is familiar with the notions of
stabilization, (ir)reducibility, and weak reducibility for Heegaard
splittings, see Scharlemann \cite{scharlemannSurvey} for the basics. For knot
manifolds, it is also worth identifying splittings that are {\it
$\gamma$--primitive}.  They possess a ``strong'' $(A_\gamma,D)$
pair: $A_\gamma$ is a vertical annulus in the compression body with
slope $\gamma$ on $\partial X$, and $D$ a disk in the handlebody so
that $|\partial A_\gamma \cap \partial D|$ = 1.   A step down from
this are splittings that are weakly $\gamma$--primitive. They
possess a ``weak'' $(A_\gamma,D)$ pair: a vertical annulus in the
compression body with slope $\gamma$ and an essential disk in the
handlebody that are disjoint. We will introduce several other
triviality conditions in later sections, identifying splittings that
are {\it padded, parallel stabilizations}, and {\it boundary
stabilizations}.

\subsection{Heegaard structure -- the Heegaard  tree and canopy}
\label{secStabTree}

What is meant by the Heegaard structure of a manifold $M$? Of
course, we  include the Heegaard genus of $M$ as well as the set of
irreducible Heegaard surfaces for $M$. A bit more general is the
{\it Heegaard tree for $M$}, $\mathcal{HT}_M$.   A vertex in the
tree, $v_\Sigma$, is the class of surfaces in $M$ isotopic to a
Heegaard surface $\Sigma$. A directed edge will point from
$v_\Sigma$ to $v_\Sigma'$ if $\Sigma'$ can be obtained from $\Sigma$
by a single stabilization.

The word tree is a bit of a misnomer.  For a manifold with fewer
than two boundary components $\mathcal{HT}_M$ is indeed a tree. But,
a Heegaard splitting $M = V \cup_\Sigma W$ induces a partition $\{
\partial M \cap V \parallel \partial M \cap W \}$ of the boundary
components.  And isotopy and stabilization, hence destabilization,
cannot change this partition.  So $\mathcal{HT}_M$ actually consists
of $2^{|\partial M| -1}$ (infinite) components, each a tree as
proved by the the Reidemeister-Singer Theorem. The leaves of the
tree are precisely the non-stabilized splittings.

We also define a somewhat finer variation, {\it the oriented
Heegaard tree for $M$}, $\mathcal{HT}_M^\pm$, where the vertices are
now defined to be isotopy classes of \emph{oriented} Heegaard
surfaces.   In this case we will assume that the compression bodies
are ordered 1 and 2 and the orientation on the Heegaard surface
points towards the second compression body. In other words, we
differentiate between the Heegaard splittings $V \cup_\Sigma W$ and
$W \cup_\Sigma V$.   Note that flipping the orientation of the
Heegaard surface swaps the partition of boundary components.   It
follows that  $\mathcal{HT}_M^\pm$ will consist of $2^{|\partial
M|}$ connected trees, double that of $\mathcal{HT}_M$.  For a knot
manifold, we can unambiguously define the {\it sign} ($\pm$) of an
oriented Heegaard surface:  A positive ($+$) Heegaard surface will
have its orientation pointing into the handlebody, and a negative
($-$) Heegaard surface will have its orientation pointing into the
compression body.

 While this may seem somewhat obvious,
it does underline an important difference between the Heegaard
structure of a knot manifold and those obtained by Dehn filling. A
knot manifold $X$ has connected boundary, so $\mathcal{HT}_X^\pm$
consists of two components, $\mathcal{HT}_X^+$ and
$\mathcal{HT}_X^-$, each homeomorphic to the tree $\mathcal{HT}_X$.
But, for a filled manifold $X(\alpha)$, the tree
$\mathcal{HT}^\pm_{X(\alpha)}$ is connected.  We will say that a
surface {\it flips} if it is isotopic to itself with reverse
orientation.  This is equivalent to an isotopy that takes the
handlebody $V$ to the handlebody $W$ for the splitting $V
\cup_\Sigma W$. In a closed manifold, there is always a Heegaard
surface that flips:  Since the oriented Heegaard tree is connected,
an oriented surface and its reverse have a common stabilization;
this surface flips.  It is an easy exercise to show that the common
stabilization of a Heegaard surface of genus $g$ and its reverse has
genus at most $2g$; hence the smallest genus of surfaces that flip
is at most twice the genus of the manifold.   But a knot manifold
never has a surface that flips, because a handlebody is never
isotopic to a compression body. In other words, for \emph{every}
filled manifold $X(\alpha)$ there are Heegaard surfaces for $X$ that
{\it flip} in $X(\alpha)$ but not in $X$.

Since we can stabilize a given splitting any number of times, each
tree defined above is infinite.   Instead of drawing
$\mathcal{HT}_X^\pm$ (upside down!), we will instead draw its {\it
canopy}, that is the smallest subset of $\mathcal{HT}_X^\pm$ that
has the same number of components as $\mathcal{HT}_X^\pm$ and that
contains all of its leaves (non-stabilized splittings). A result of
Li \cite{li} shows that the stabilization tree has an infinite
canopy only if the manifold contains a closed essential surface.
Examples of canopies for $\mathcal{HT}_M^\pm$ are drawn in \fullref{secTorusKnotFillings}.

\section{Dehn Filling on Torus Knots}
\label{secTorusKnotFillings}

In this section we will review the Heegaard structure of torus knot
exteriors and the manifolds that can be obtained by Dehn filling on
them. Fortunately there has been a lot of work done in this area,
and we know the Heegaard tree $\mathcal{HT}_{M}^\pm$ up to isotopy
for all torus knot exteriors and almost all manifolds that can be
obtained by Dehn filling on a torus knot exterior.   The sole
exception is a restricted class of connected sums of lens spaces
whose Heegaard structure is known up to homeomorphism but not up to
isotopy.  These are discussed in \fullref{secConnectedSum}.

First, we fix notation.  Let $T$ be a Heegaard torus in $S^3$, it
separates $S^3$ into two solid tori that we will denote by $V_i$ and
$V_o$. Let $\mu$ and $\lambda$ denote the meridians of $V_i$ and
$V_o$, respectively.  Then the curve $T_{p,q} = p\lambda + q \mu$ is
a $(p,q)$--torus knot in $S^3$. The exterior, $X = S^3 - N(T_{p,q})$
is a Seifert fibered space over the disk with two exceptional fibers
$f_i$ and $f_o$.  A regular neighborhood of $f_i$ is the solid torus
$V_i$ with a $(p,q)$--fibering by regular fibers and a regular
neighborhood of $f_o$ is the solid torus $V_o$ with a $(q,p)$
fibering by regular fibers (see Scott \cite{scottGeometries} and Jaco
\cite{jacoBook}).

\subsection{Our examples:}
\label{secOurExamples}

 We will restrict our attention to $(p,q)$--torus knot
exteriors that satisfy two conditions:

\begin{enumerate}
\item $p \not \equiv \pm 1 \pmod q$ and $q \not \equiv
\pm 1 \pmod p$
\item $q^2 \not \equiv \pm 1 \pmod p$ and $p^2 \not \equiv \pm 1
\pmod q$
\end{enumerate}

The first condition rules out torus knot exteriors with fewer than
three non-isotopic tunnels.  This restriction keeps our listing of
$\mathcal{HT}_M^\pm$'s a bit shorter, but the excluded knots can be
analyzed in the same manner. The second rules out fillings that
produce a connected sum of lens spaces whose Heegaard structure is
known only up to homeomorphism. As mentioned above, this will be
discussed further in \fullref{secConnectedSum}.

\subsection{Heegaard structure of {$\mathbf{\{\text{pair of pants}\} \times
S^1}$}} \label{secPOP}

\begin{figure}[ht!]
\begin{center}
\includegraphics[height=1.25in]{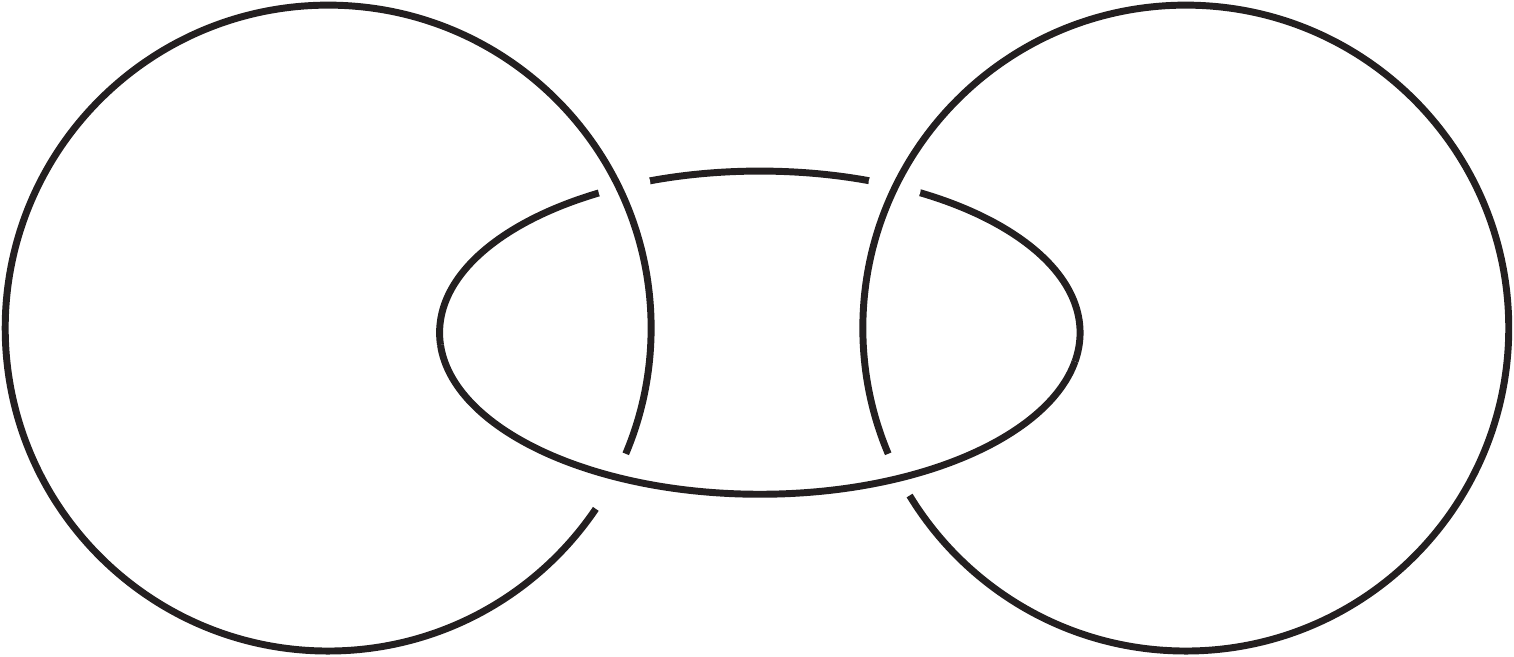} \qquad
\labellist\small
\pinlabel {$1$} at 200 480
\pinlabel {$2$} at 238 378
\pinlabel {$3$} at 355 378
\endlabellist
\includegraphics[height=1.65in]{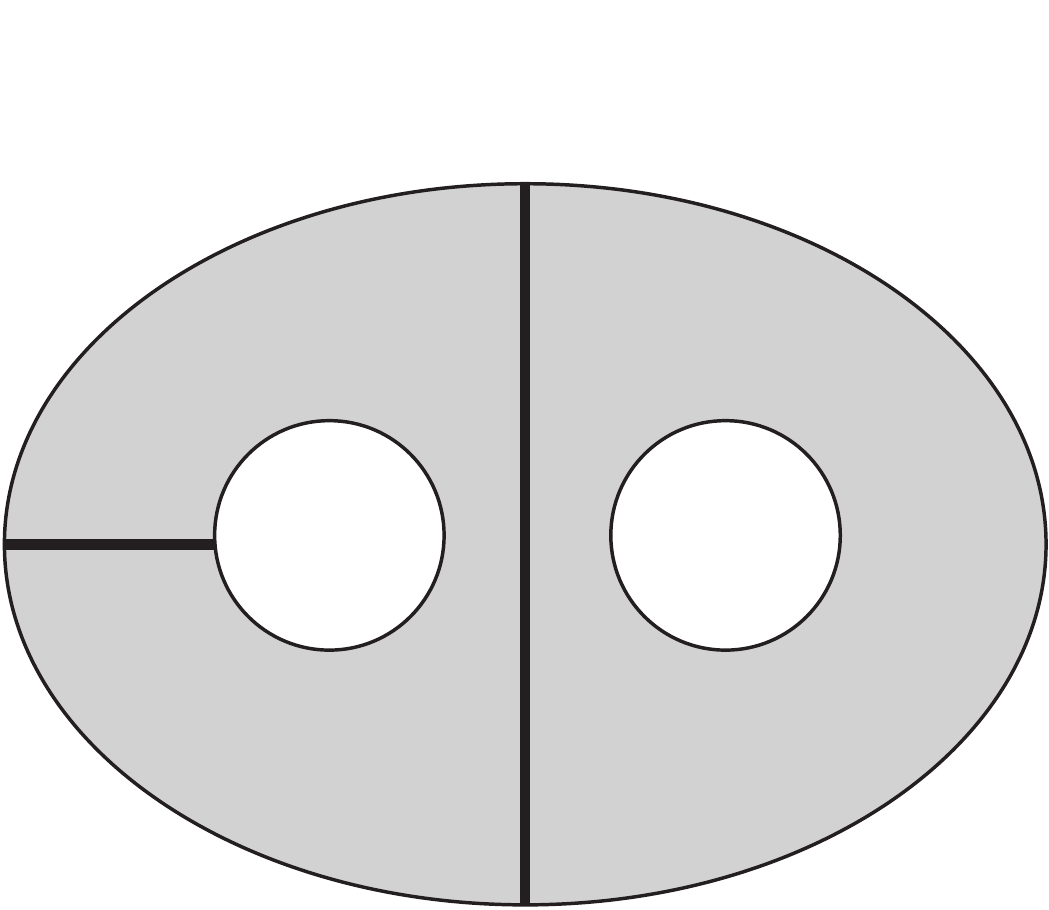}

\caption{The link manifold $X = S^3 - N(\text{chain with 3
components})$ is homeomorphic to $\text{\{pair of pants\}} \times
S^1$. The arcs $a_{12}$ and $a_{11}$ in the pair of pants $P$.}
\label{figChain}
\end{center}
\end{figure}

Heegaard splittings of Seifert fibered spaces with boundary are
vertical (see Schultens \cite{schultensSFSBoundary}).   A {\it vertical
splitting} is a Heegaard splitting for the Seifert fibered space
that is also a Heegaard surface for the product manifold obtained by
drilling out all of the exceptional fibers.   To understand
splittings of the torus knot exterior, we look to the corresponding
product manifold $P \times S^1$, where $P$ is a pair of pants.  This
manifold is also homeomorphic to the exterior of the three component
chain in $S^3$, pictured in \fullref{figChain}.  We have already
noted that Heegaard splittings of a manifold with boundary induce
partitions of the boundary components. Heegaard splittings of $P
\times S^1$ are special, because any partition identifies, up to
isotopy, a unique irreducible splitting (see \cite{schultensSFSBoundary}).

While the product structure on $P \times S^1$ is not unique, the
Seifert fibering is unique (see Jaco \cite{jacoBook}).   Let $A_{ij}$ denote
the unique essential annulus joining $\partial X_i$ and $\partial
X_j$, and $a_{ij} = A_{ij} \cap P$ a spanning arc for $A_{ij}$.  A
different choice of $P$ will yield a different, but handle slide
equivalent, spanning arc $a_{ij}$.  Any curve in $\partial X_i$
which meets a regular fiber once will be referred to as a {\it dual
curve}.   Dual curves determine the possible boundary slopes for
different choices for $P$.

There are three genus two splittings (six when oriented), each
identified by a partition of boundary components: $\Sigma_{12}^+
\leftrightarrow \{\partial X_1,\partial X_2
\parallel
\partial X_3\}, \Sigma_{13}^+ \leftrightarrow\{\partial X_1,\partial X_3
\parallel
\partial X_2\}$ and $\Sigma_{23}^+ \leftrightarrow \{\partial X_2,\partial X_3
\parallel
\partial X_1\}$.   The compression bodies corresponding to
$\Sigma_{ij} \leftrightarrow \{\partial X_i, \partial X_j \parallel
\partial X_k\}$ are isotopic to $N(\partial X_i \cup a_{ij} \cup
\partial X_j)$ and $N(\partial X_k \cup a_{kk} $). There is also an irreducible
splitting of genus three identified by the partition $\Sigma_{123}^+
\leftrightarrow \{\partial X_1,
\partial X_2, \partial X_3 \parallel \emptyset\}$.  In that case the compression body
is given by  $N(\partial X_1 \cup a_{12} \cup \partial X_2 \cup
a_{23} \cup
\partial X_3 $).

The canopy of $\mathcal H_{P \times S^1}^\pm$ is indicated in \fullref{figPOP}.  It has exactly one non-stabilized oriented Heegaard
surface for each oriented partition of boundary components.

\begin{figure}[ht!]
$$\begin{array}{ccccccccccc}
g= 3  &&& ~\Sigma_{123}^+~&&&&&&&~\Sigma_{123}^-~\\[2ex]
g=2&&&& ~\Sigma_{12}^+~&~\Sigma_{13}^+~& ~\Sigma_{23}^+ ~&
~\Sigma_{23}^-~&~\Sigma_{13}^-~& ~\Sigma_{12}^- ~
\end{array}$$

 \caption{Canopy of $\mathcal{H_X^+}$ for  $\{\text{pair of
 pants}\} \times S^1$.   Up to isotopy, there is a single non-stabilized oriented Heegaard surface for each
 ordered partition of boundary components, for example $\Sigma_{12}^+ \leftrightarrow \{\partial X_1, \partial X_2 \parallel \partial
 X_3\}$.}
 \label{figPOP}
\end{figure}

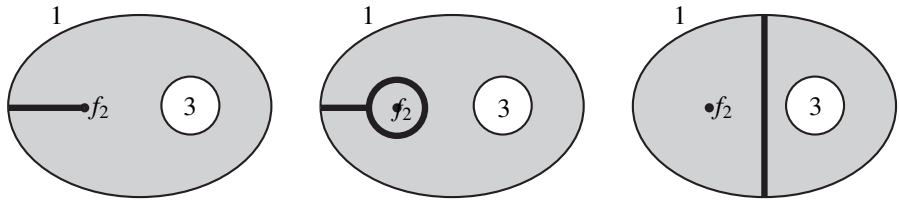
\begin{figure}[ht!]
\begin{center}
\labellist\small
\pinlabel {$1$} at 109 425
\pinlabel {$f_2$} [l] at 122 381
\pinlabel {$3$} at 172 382
\pinlabel {$1$} at 256 425
\pinlabel {$f_2$} [l] at 265 380
\pinlabel {$3$} at 320 381
\pinlabel {$1$} at 403 425
\pinlabel {$f_2$} [l] at 416 381
\pinlabel {$3$} at 467 382
\endlabellist
\includegraphics[height=1in]{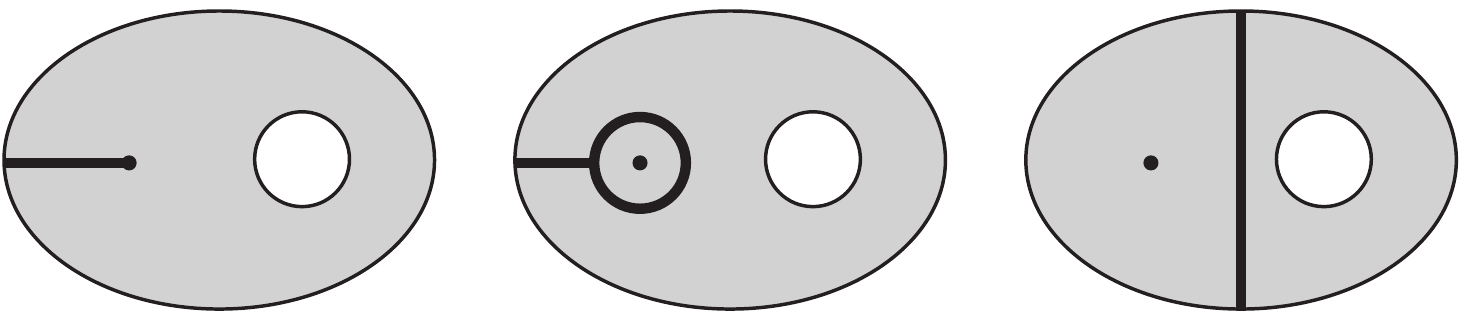}

\caption{``Flipping'' the fiber $f_2$: from $a_{12} \cup f_2$ to
{\it eyehook} to $a_{11}$. Requires that $b \equiv \pm 1 \pmod a$
where $(a,b)$ are the Seifert invariants of the fiber being flipped.}
\label{figFlipFiber}
\end{center}
\end{figure}

Of course, each of these splitting will also be a splitting after we
fill in any or all of the exceptional fibers.   In that case the
splitting will be identified by partitions of boundary components
{\it and} exceptional fibers $f_i$, for example $\{f_i,\partial X_j
\parallel
\partial X_k\}$.  After filling, it is possible that splittings
corresponding to distinct  partitions now become isotopic. Fill
$\partial X_2$ and consider the Heegaard surface $\Sigma_{12}^+$
inducing the partition $\{ \partial X_1, f_2
\parallel
\partial X_3 \}$.   In this case the first compression body is a
regular neighborhood of $\partial X_1 \cup a_{12} \cup f_2$. Suppose
that the Seifert invariants of the fiber are $(p,q)$ and that $q
\equiv \pm 1 \pmod p$. This implies that we can find some {\it
longitude}, a curve meeting the meridional disk of the attached
solid torus once, that is also a dual curve meeting a regular fiber
exactly once. In other words, we can isotope $a_{12} \cup f_2$ to
appear as an eyehook on $P$, as in the central picture of \fullref{figFlipFiber}. Sliding the foot of the circle to $\partial X_1$
does not change the isotopy class of the Heegaard surface and
demonstrates that this splitting is equivalent to the splitting
induced by $a_{11}$.   We have changed the partition from $\{
\partial X_1, f_2 \parallel
\partial X_3 \}$ to $\{\partial X_1
\parallel f_2,\partial X_3 \}$, demonstrating an isotopy between $\Sigma_{12}^+$ and $\Sigma_{23}^-$.
The isotopy {\it flips} $f_2$, moving it from one side of the
partition to the other.    Any fiber with Seifert invariants $(p,q)$
so that $q \equiv \pm 1 \pmod p$ can be flipped. In essence, such
fibers should be left out of the partition altogether.   In fact,
this is the only way that vertical splittings of a Seifert fibered
spaces become isotopic (see Schultens \cite[Theorem~5.1]{schultensSFSBoundary}).

The genus three splitting is very fragile.   While irreducible, it
is boundary stabilized (see Moriah \cite{moriahWR}).  The notion
of boundary stabilization will be discussed in greater detail in
\fullref{secDestab}. In fact, it can be viewed as a boundary stabilization
of each of the genus two splittings. This implies that it will destabilize
after any filling on any one of the three boundary components.

\subsection{Heegaard
structure of the $\mathbf{(p,q)}$--torus knot exterior}

Genus two Heegaard splittings of torus knot exteriors were
originally classified by the first author \cite{moriahSFS} (see also
Boileau, Rost and Zieschang \cite{brz}).
Since the torus knot exterior is a Seifert fibered space with
boundary, any irreducible Heegaard splitting is vertical, hence
isotopic to one of the three (unoriented) irreducible genus two
Heegaard splittings  of $\{\textit{pair of pants}\} \times S^1$,
discussed in the previous section. We have pictured two of these as
tunnels for the torus knot in \fullref{figTunnels}, they are:

\begin{figure}[ht!]
\begin{center}
\includegraphics[height=1.75in]{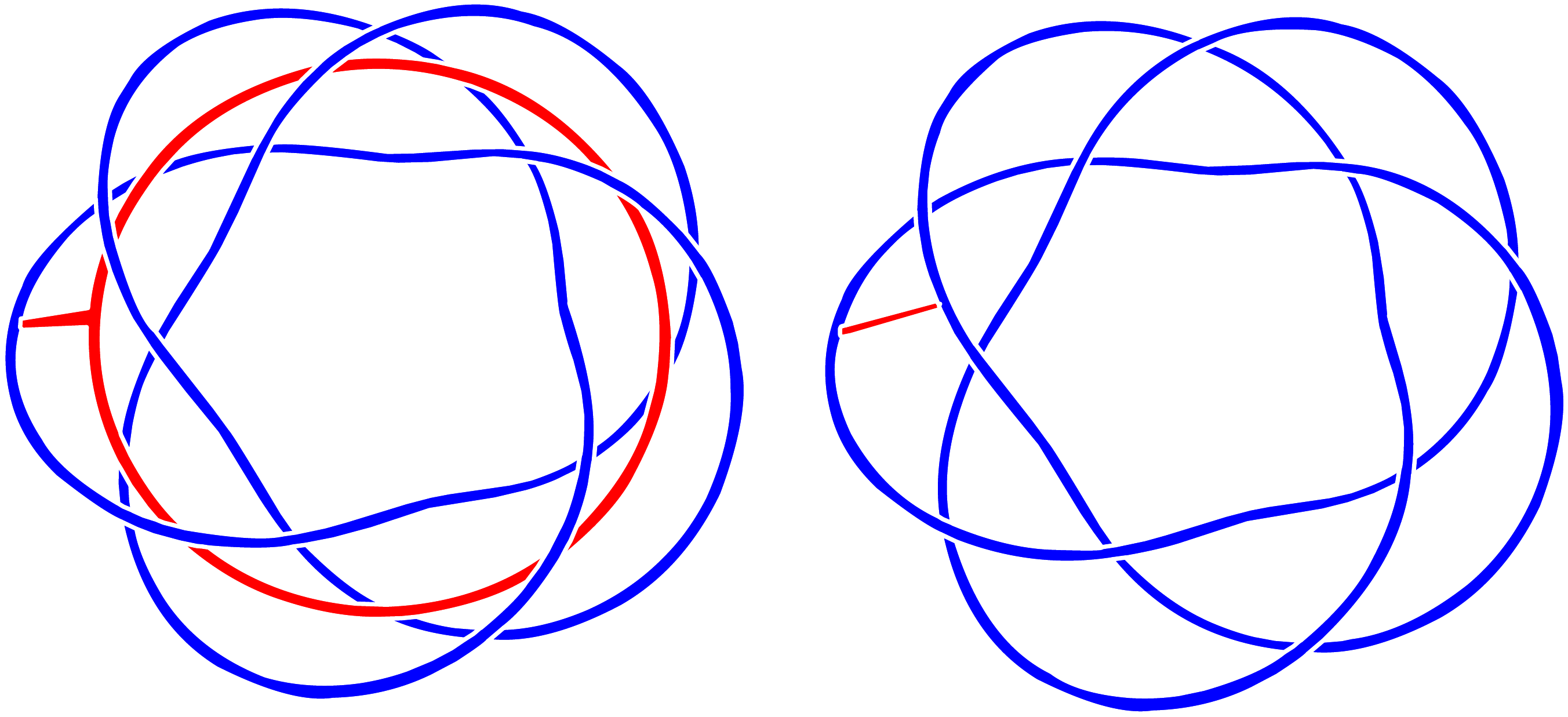}

\caption{The inner and middle tunnels for a torus knot}
\label{figTunnels}

\end{center}

\end{figure}

\begin{enumerate}
\item The inner tunnel -- the inner exceptional fiber $f_i$ joined to the knot via a vertical arc (a spanning arc for the
annulus between $\partial X$ and $f_i$) - $\Sigma_i^+
\leftrightarrow \{ f_i,
\partial X
\parallel f_o\}$.
\item The outer tunnel -- the outer exceptional fiber $f_o$ joined to the knot via a vertical arc (a spanning arc for the
annulus running between $\partial X$ and $f_o$) - $\Sigma_o^+
\leftrightarrow \{ f_o, \partial X
\parallel f_i\}$.
\item The middle tunnel -- a spanning arc for the {\it cabling annulus} $A = T -
N(T_{p,q})$- $\Sigma_m^+ \leftrightarrow \{ \partial X
\parallel f_i, f_o\}$.
\end{enumerate}

These three splittings are distinct up to isotopy (and
homeomorphism), unless $|p-q|=1$, in which case all three are
isotopic; or $|p-q| \neq 1$, but $p \equiv \pm 1\mod q$ or $q \equiv
\pm 1 \mod p$, in which case the middle splitting is isotopic to the
inner or outer splitting, respectively (see Boileau, Rost and Zieschang
\cite{brz} or Moriah \cite{moriahSFS}, and the previous section).  All (oriented) Heegaard
surfaces of closed Seifert fibered spaces are equivalent after one
stabilization (see Schultens \cite{schultensStabilization}), therefore
the canopy of the Heegaard
tree is:

\begin{figure}[ht!]
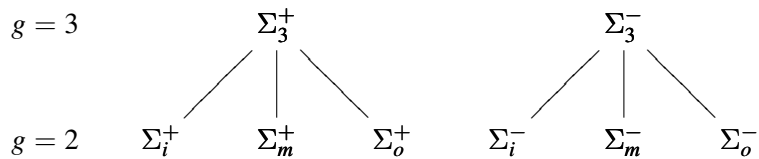

$$\bfig
  \morphism//<0,-400>[g=3`g=2;]
  \Atrianglepair(400,-400)/-`-`-``/<400,400>[\Sigma_3^+`
    \Sigma_i^+`\Sigma_m^+`\Sigma_o^+;````]
  \Atrianglepair(1600,-400)/-`-`-``/<400,400>[\Sigma_3^-`
    \Sigma_i^-`\Sigma_m^-`\Sigma_o^-;````]
  \efig$$
 \caption{Canopy of $\mathcal{H_X^+}$ for a
$(p,q)$--torus knot exterior $X$. (Unless $|p-q|=1$, in which case
$\Sigma_i^+ = \Sigma_m^+ = \Sigma_o^+$; or $|p-q| \neq 1$, but $p
\equiv \pm 1\mod q$ or $q \equiv \pm 1 \mod p$, in which case
$\Sigma_i^+ = \Sigma_m^+$ or $\Sigma_o^+ = \Sigma_m^+$,
respectively). }
\end{figure}

\subsection{The manifolds obtained by filling}

The manifolds obtained by a Dehn filling on a torus knot exterior
were classified by Moser \cite{moser}, whose theorem is rephrased
slightly here:

\begin{thm}[Moser]
\label{thmMoser} Suppose that $\frac{r}{s}$--Dehn filling is
performed on a non-trivial $(p,q)$--torus knot and let $a =
\Delta_A(\frac{pq}{1},\frac{r}{s}) = pqs-r$ be the algebraic
intersection number between the slope of a regular fiber and the
meridian of the attached solid torus.   The type of the filled
manifold $X(\frac{r}{s})$ depends on $a$:
\begin{enumerate}
\item $a=0 ~\implies X(\frac{r}{s}) = L(p,q) \# L(q,p)$, a connected sum of lens
spaces.
\item $|a|=1 \implies X(\frac{r}{s}) = L(|r|,sq^2)$, a lens space.
\item $|a|>1 \implies X(\frac{r}{s}) = \textit{SFS}\{S^2|(p,q),(q,p),(a,b)\}$, a Seifert fibered space over $S^2$ with three exceptional
fibers with Seifert invariants  $(p,q)$, $(q,p)$ and $(a,b)$.
\end{enumerate}
The Seifert invariants are not normalized and $b =
\Delta_A(\frac{pq}{1},\frac{t}{u})=pqu-t$, where $\frac{t}{u}$ is
the slope of any longitude for the attached solid torus.
\end{thm}

{\rm
\begin{rmk}
\rm Note that $a=0 \Leftrightarrow \frac{r}{s} = \frac{pq}{1}$,
$|a|=1 \Leftrightarrow \frac{r}{s} \in L_{pq/1},$ and $|a|>1$
otherwise.
\end{rmk}
}

\subsection{Heegaard structure of $\mathbf{S^3}$}
Of course, $S^3$ appears as a special case of \fullref{thmMoser}
(2), the $\frac{1}{0}$--filling produces $L(1,0)$, that is, $S^3$.
Waldhausen \cite{waldhausen:1967} showed $S^3$ has a unique
non-stabilized Heegaard surface, an $S^2$ splitting $S^3$ into two
balls. It is also easy to see that this surface flips. Each ball is
a regular neighborhood of a point in its interior, and any two
neighborhoods of points are isotopic in $S^3$.  This defines an
isotopy that reverses the orientation of the Heegaard sphere.
Therefore the canopy of the oriented Heegaard tree for $S^3$ is as
simple as possible: It is a single point.

\begin{figure}[ht!]
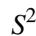

$$S^2$$
 \caption{Canopy of $\mathcal{HT}_{S^3}^\pm$ is a single point, $\alpha = \frac{1}{0}$.}
\end{figure}

\subsection{Heegaard structure of $\mathbf{L(|r|,sq^2)}$}

A lens space also has a unique non--stabilized Heegaard surface, a
torus $T$ (see Bonahon and Otal \cite{bonahon-otal}). The splitting does not flip
unless the cores of the solid tori are isotopic, hence equivalent as
generators of the fundamental group.   This will occur precisely for
those lens spaces homeomorphic to $L(k,1)$ for some $k \in
\mathbb{N}$.  A check will reveal that such a manifold can not be
obtained from surgery on a torus knot, except for the
$\frac{5}{1}$--filling on the $(3,2)$--torus knot that which
produces $L(5,4)$. So for any other lens space, including all of
those we are considering, we have:

\begin{figure}[ht!]
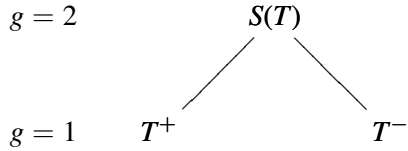

$$\bfig
  \morphism//<0,-400>[g=2`g=1;]
  \Atriangle(400,-400)/-`-`/<400,400>[S(T)`T^+`T^-;``]
  \efig$$
\caption{Canopy of $\mathcal{HT}_{L(p,q)}^\pm, q \not \equiv \pm 1 \pmod
p$, $\alpha \in L_{pq/1} - \{1/0\}$.}
\end{figure}

\subsection{Heegaard structure of $\mathbf{L(p,q)\#L(q,p)}$}
\label{secConnectedSum}

The Haken Lemma implies that any Heegaard splitting of
$L(p,q)\#L(q,p)$ is the connected sum of a Heegaard torus $T_1$ for
$L(p,q)$ and a Heegaard torus $T_2$ for $L(q,p)$.   For the
connected sum, we can form two non--oriented Heegaard surfaces,
$T_1^+ \# T_2^+$ and $T_1^+ \# T_2^-$, or four oriented Heegaard
surfaces, since $(T_1^+ \# T_2^+)^- = (T_1^- \# T_2^-)$ and $(T_1^+
\# T_2^-)^- = T_1^- \# T_2^+$. Of course, it is possible that two or
more of these surfaces are isotopic, and that will definitely be the
case when a Heegaard torus for one of the summand flips, that is, when
$q \equiv \pm 1 \pmod p$ or $p \equiv \pm 1 \pmod q$. However,
Engmann \cite{engmann} has shown that the connected sum of lens
spaces $L(p_1,q_1) \# L(p_2,q_2)$ has four distinct oriented
Heegaard surfaces up to homeomorphism, hence isotopy, unless: a)
$q_i^2 \equiv \pm 1 \pmod {p_i}$ for some $i$, or b) $p_1=p_2$ and
$q_1q_2^{-1} \equiv \pm 1 \pmod {p_1}$. Since $p \neq q$, we are
concerned only with a). In that case, it is clear that the the
Heegaard surfaces $T_i^+$ and $T_i^-$ for $L(p_i,q_i)$ are
homeomorphic, reducing the number of non-homeomorphic splittings for
the connected sum.  However, unless some $q_i \equiv \pm 1 \pmod
{p_i}$, it is not clear that these splittings are isotopic and we
conjecture that they are not.   For torus knots satisfying the
conditions 1) and 2) in \fullref{secOurExamples},
$\mathcal{HT}_{X(\alpha)}^\pm, \alpha = \frac{pq}{1}$ is indicated
in the following figure.   We also conjecture that this represents
$\mathcal{HT}_{X(\alpha)}^\pm, \alpha = \frac{pq}{1}$ assuming only
condition (1).

\begin{figure}[ht!]
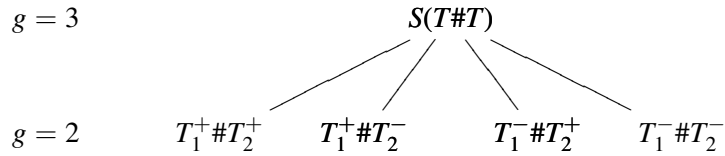

$$\bfig
  \morphism//<0,-400>[g=3`g=2;]
  \morphism(600,-400)/-/<800,400>[T_1^+{\#}T_2^+`\phantom{S(T{\#}T)};]
  \Atriangle(1100,-400)/-`-`/<300,400>[S(T{\#}T)`
    T_1^+{\#}T_2^-`T_1^-{\#}T_2^+;``]
  \morphism(2200,-400)/-/<-800,400>[T_1^-{\#}T_2^-`\phantom{S(T{\#}T)};]
  \efig$$
\caption{Canopy for $\mathcal{HT}_{L(p,q)\#L(q,p)}^\pm$, obtained by
  filling along $\alpha = \frac{pq}{1}$.  This holds when $q^2 \not
  \equiv \pm 1 \pmod p, p^2 \not \equiv \pm 1 \pmod q$ and we conjecture
  that it holds with the weaker assumption $q \not \equiv \pm 1 \pmod p,
  p \not \equiv \pm 1 \pmod q$.}

 \label{figLpqLqp}
\end{figure}

\subsection{Heegaard structure of
$\mathbf{\textit{SFS}\{S^2~|~(p,q),(q,p),(a,b)\}}$}

Heegaard splittings of Seifert fibered spaces are either vertical or
horizontal (see Moriah--Schultens \cite{ms} and Schultens
\cite{schultensExceptional}). The vertical splittings of a Seifert
fibered spaces over $S^2$ with three exceptional fibers are inherited
from the torus knot exterior: $\Sigma_i, \Sigma_o$ and $\Sigma_m$. All
three are genus two, therefore minimal genus and irreducible. There may
or may not be an irreducible horizontal splitting, but if it exists it
is unique up to isotopy (see Moriah and Schultens \cite{ms}, Sedgwick
\cite{sedgwickSFS} and Bachman--Derby-Talbot \cite{bachmanSFS}).

As already noted, the three splittings will be distinct up to
isotopy unless $b \equiv \pm 1 \pmod {|a|}$ (see Moriah
\cite[Theorem~1]{moriahSFS} and Schultens
\cite[Theorem~5.1]{schultensSFSBoundary}). In
that case, the inner and outer splittings are isotopic. More
specifically, $\Sigma_i^+$ is isotopic to $\Sigma_o^-$  (see \fullref{secIsotopic}). This occurs when there is a curve with slope
$\frac{t}{u}$ which is a longitude for the filling solid torus and
that meets the regular fiber of the Seifert fibration once, that is,
$\Delta(\frac{t}{u},\frac{r}{s})=1$.  This is precisely the set of
slopes on the line of lines $LL_{pq/1}$. Equivalently, it is the set
of slopes with distance $2$ from the slope of the regular fiber
$\frac{pq}{1}$ in the Farey graph.

\begin{figure}[ht!]
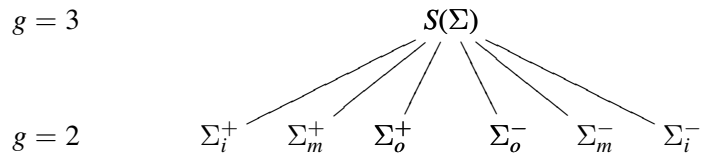

$$\bfig
  \morphism//<0,-400>[g=3`g=2;]
  \morphism(600,-400)/-/<800,400>[\Sigma_i^+`\phantom{S(\Sigma)};]
  \morphism(900,-400)/-/<500,400>[\Sigma_m^+`\phantom{S(\Sigma)};]
  \Atriangle(1200,-400)/-`-`/<200,400>[S(\Sigma)`
    \Sigma_o^+`\Sigma_o^-;``]
  \morphism(1900,-400)/-/<-500,400>[\Sigma_m^-`\phantom{S(\Sigma)};]
  \morphism(2200,-400)/-/<-800,400>[\Sigma_i^-`\phantom{S(\Sigma)};]
  \efig$$
\caption{Canopy for the ``generic'' Seifert fibered space obtained
  by filling a torus knot exterior.  There are three non-isotopic
  vertical splitting. $\mathcal{HT}_{X(\alpha)}^\pm, \alpha \not \in
  L_{0/1} \cup LL_{pq/1}$. }
\label{figSFSTypical}
\end{figure}

There is only one line, $L_{0/1}=\{\frac{1}{n}| n \in \mathbb{Z}\}$,
of fillings which produce a manifold with a strongly irreducible
{\it horizontal} splitting. Recall that the torus knot exterior
fibers (uniquely) over the circle, $X = F \tilde{\times} S^1$. In
$X$ a neighborhood of $F$ is a handlebody, as is its complement.  If
the meridional slope $\frac{r}{s}$ for the filling solid torus meets
the slope of $\partial F = \frac{0}{1}$ once, then the solid torus
can be glued to either handlebody and in both cases the result will
still be a handlebody, hence we get a Heegaard splitting of
$X(\frac{1}{n})$. It has genus $g=2g(F) = 2(p-1)(q-1)$. The
splitting is (strongly) irreducible if and only if $|a| = |pqn-1|
> \lcm(p,q) = pq$ \cite{sedgwickSFS}.  That
is, the horizontal splitting is strongly irreducible unless
$\frac{1}{n}$ is one of $\frac{1}{0}$ or $\frac{1}{1}$. (Since
$X(\frac{1}{0}) = S^3$, we already know that the horizontal
splitting reduces in that case.)

\begin{figure}[ht!]
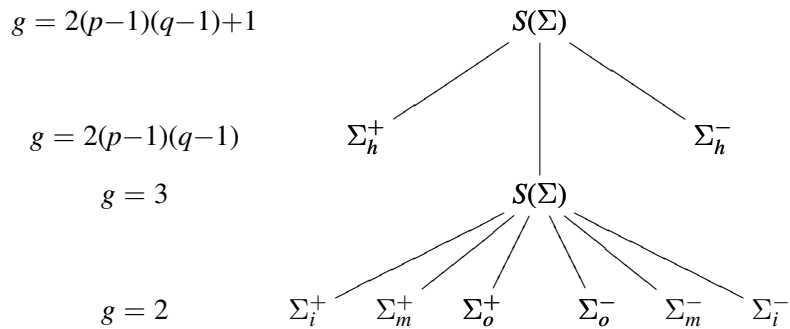

$$\bfig
  \morphism//<0,-400>[g=2(p{-}1)(q{-}1){+}1`g=2(p{-}1)(q{-}1);]
  \Atriangle(800,-400)/-`-`/<600,400>[S(\Sigma)`
    \Sigma_h^+`\Sigma_h^-;``]
  \morphism(1400,0)/-/<0,-600>[\phantom{S(\Sigma)}`\phantom{S(\Sigma)};]
  \morphism(0,-600)//<0,-400>[g=3`g=2;]
  \morphism(600,-1000)/-/<800,400>[\Sigma_i^+`\phantom{S(\Sigma)};]
  \morphism(900,-1000)/-/<500,400>[\Sigma_m^+`\phantom{S(\Sigma)};]
  \Atriangle(1200,-1000)/-`-`/<200,400>[S(\Sigma)`
    \Sigma_o^+`\Sigma_o^-;``]
  \morphism(1900,-1000)/-/<-500,400>[\Sigma_m^-`\phantom{S(\Sigma)};]
  \morphism(2200,-1000)/-/<-800,400>[\Sigma_i^-`\phantom{S(\Sigma)};]
  \efig$$
 \caption{Canopy for Seifert fibered spaces possessing three
 vertical splittings and an irreducible horizontal splitting.
 $\mathcal{HT}_{X(\alpha)}^\pm, \alpha \in L_{0/1} -
 \{\frac{1}{0},
 \frac{1}{1},\frac{-1}{1}\}$.}
 \label{figSFSwHorizontal}
\end{figure}

Although there are Seifert fibered spaces which have infinitely many
non-isotopic horizontal splittings, they are always obtained by Dehn
twisting in non-separating tori \cite{bachmanSFS}. However,
the Seifert fibered spaces in question here have no essential tori,
so the fiber $F$, hence the horizontal splitting, are both unique up
to isotopy.

As we have already noted, Schultens has shown that Heegaard splittings
of Seifert fibered spaces are equivalent after one stabilization
\cite{schultensStabilization}. Her argument applies to oriented Heegaard
surfaces of closed Seifert fibered spaces.

\begin{figure}[ht!]
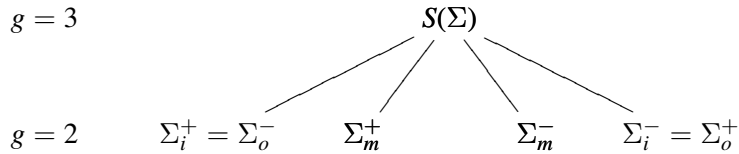


$$\bfig
  \morphism//<0,-400>[g=3`g=2;]
  \morphism(600,-400)/-/<800,400>[\Sigma_i^+=\Sigma_o^-`\phantom{S(\Sigma)};]
  \Atriangle(1100,-400)/-`-`/<300,400>[S(\Sigma)`\Sigma_m^+`\Sigma_m^-;``]
  \morphism(2200,-400)/-/<-800,400>[\Sigma_i^-=\Sigma_o^+`\phantom{S(\Sigma)};]
  \efig$$
    \caption{Canopy for the Seifert fibered space with two
    non-isotopic vertical splittings and no horizontal splitting,
    $\mathcal{HT}_{X(\alpha)}^\pm, \alpha \in  LL_{pq/1} - \{
\frac{-1}{1}, \frac{1}{0}, \frac{1}{1}\}$.}
\label{figSFS2Vertical}
\end{figure}

\begin{figure}[ht!]
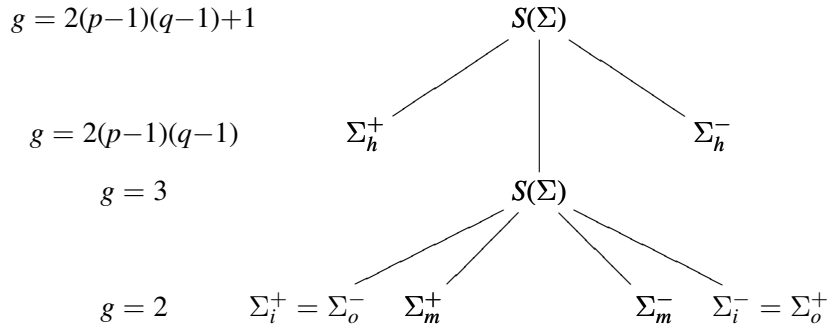


$$\bfig
  \morphism//<0,-400>[g=2(p{-}1)(q{-}1){+}1`g=2(p{-}1)(q{-}1);]
  \Atriangle(800,-400)/-`-`/<600,400>[S(\Sigma)`
    \Sigma_h^+`\Sigma_h^-;``]
  \morphism(1400,0)/-/<0,-600>[\phantom{S(\Sigma)}`\phantom{S(\Sigma)};]
  \morphism(0,-600)//<0,-400>[g=3`g=2;]
  \morphism(600,-1000)/-/<800,400>[\Sigma_i^+=\Sigma_o^-`\phantom{S(\Sigma)};]
  \Atriangle(1000,-1000)/-`-`/<400,400>[S(\Sigma)`
    \Sigma_m^+`\Sigma_m^-;``]
  \morphism(2200,-1000)/-/<-800,400>[\Sigma_i^-=\Sigma_o^+`\phantom{S(\Sigma)};]
  \efig$$
\caption{Canopy for Seifert fibered space with two non-isotopic vertical
  splittings and an a horizontal splitting.  $\mathcal{HT}_{X(\alpha)}^\pm,
  \alpha = \frac{-1}{1}$.}
\label{figSFSHorizontalAnd2Vertical}
\end{figure}

It follows that we have four possibilities for the stabilization tree
$\mathcal H_{X(\alpha)}^\pm$ when the filled manifold is a Seifert fibered
space that is not a lens space or $S^3$.   In the ``generic'' case,
there will be three non-isotopic vertical splittings (six oriented).
Exceptions will occur when the filled manifold has only two vertical
surfaces up to isotopy ($\frac{r}{s} \in LL_{pq/1}$), the filled manifold
contains an irreducible horizontal splitting of higher genus ($\frac{r}{s}
\in L_{0/1}$), or both ($\frac{r}{s} \in LL_{pq/1} \cap L_{0/1}$).
We leave the fact that $LL_{pq/1} \cap L_{0/1} = \{ \frac{-1}{1},
\frac{1}{0}, \frac{1}{1}\} $ as an exercise. Note however, that the
$\frac{1}{0}$--filling produces $S^3$ and the horizontal splitting is
reducible after the $\frac{1}{1}$--filling, so the $\frac{-1}{1}$--filling
is the unique (!) slope yielding two vertical splitting and a
horizontal splitting.  The four possible canopies are indicated in
Figures \ref{figSFSTypical}--\ref{figSFSHorizontalAnd2Vertical}.

\section{New Heegaard surfaces -- the framework}
\label{secNewSplittings}

Suppose that we are given a surface $\Sigma \subset X$ and are asked
for the set of filled manifolds for which it is a Heegaard surface.
First, if $\Sigma$ is a Heegaard surface for any $X(\alpha)$, then
$\Sigma$ must separate $X$ into two components $V$ and $W$, where
$W$ is a handlebody and $V$ is a {\it punctured handlebody}, a
handlebody with the neighborhood of a knot removed. Then $\Sigma$ is
a Heegaard surface for $X(\alpha)$ if and only if $V(\alpha)$ is a
handlebody.

We can use the results of Culler-Gordon-Luecke-Shalen \cite{cgls} and
Wu \cite{wu} to address this situation.  To implement their results
it is required that $\partial V - \partial X$ is incompressible, and
we have not assumed that. We can reduce to that situation by first
maximally compressing $\Sigma$ in $V$ to obtain $\Sigma' \subset V$ an
incompressible surface bounding a punctured handlebody $V' \subset V$.
If $\Sigma'$ is peripheral, then $\Sigma$ was itself a Heegaard surface
for the knot exterior and thus for every filling on $X$.

Suppose that $\Sigma'$ is not peripheral and that $V'$ contains an
incompressible annulus $A_\sigma$ with one boundary component in
$\Sigma$ and the other a curve of slope $\sigma$ in $\partial X$. If
$\Sigma$ bounds a handlebody in $X(\alpha)$ for $\alpha \neq
\sigma$, then $\Sigma'$ compresses in $X(\alpha)$ and by
\cite[Theorem~2.4.3]{cgls}, $\alpha \in L_\sigma$ and $\Sigma$ is a
Heegaard surface for every $\alpha \in L_\sigma$. In this case we
will call $\Sigma$ a {\it horizontal Heegaard surface}.  It is a
Heegaard surface for every filling $X(\alpha), \alpha \in L_\sigma$.
By sliding the annulus off the scars from the disk compressions and
reattaching the disks, we can assume the annulus runs between
$\Sigma$ and $\partial X$.  For any slope $\alpha \in L_\sigma$,
this annulus defines an isotopy of the core of the attached solid
torus into the surface $\Sigma$.

If there is no incompressible annulus joining $\Sigma$ and $\partial
X$ and $\Sigma$ is a Heegaard surface for $X(\alpha_1)$ and
$X(\alpha_2)$, then the results of Wu \cite{wu} imply that
$\Delta(\alpha_1,\alpha_2) = 1$.  There are therefore at most three
fillings for which $\Sigma$ is a Heegaard surface.

This analysis does not solve our problem because we do not have a
nice short list of candidate $\Sigma$'s to become Heegaard surfaces.
It does, however, provide a useful framework that is used in various
papers on this subject. Let $\Sigma$ be a Heegaard surface for a
manifold $X(\alpha)$ that is obtained by performing a Dehn filling
on a knot manifold $X$. We regard the core curve $\gamma$ of the
attached solid torus as a knot in $X(\alpha)$.  The core of the
attached solid torus $\gamma$ and the Heegaard surface $\Sigma$ can
have one of three possible relationships after performing isotopies
in $X(\alpha)$:

\begin{description}
\item[(C)ore] $\gamma$ is isotopic into $\Sigma$ and can be further isotoped so that $\Sigma$ is a Heegaard
surface for the knot exterior $X$.
\item[(H)orizontal but not a core] $\gamma$ is isotopic into
$\Sigma$ but cannot be isotoped so that $\Sigma$ is a Heegaard
surface for $X$.
\item[(N)ot level] $\gamma$ cannot be isotoped into $\Sigma$.
\end{description}

Case C describes an old splitting not a new one. For a new splitting
$\Sigma$, the issue is then whether or not the core of the attached
solid torus $\gamma$ is isotopic into $\Sigma$, case H, or not, case
N. We would like to limit the set of slopes for which condition H or
N occurs.  Condition N is also referred to as a ``bad'' filling in
Rieck--Sedgwick \cite{rs1}.

Which of the Heegaard surfaces from the previous section are {\it
new}, that is, not isotopic in the filled manifold to a Heegaard
surface for the torus knot exterior?   By definition, the vertical
splittings are not new as they are also splittings of the knot
exterior. Every other non-stabilized splitting is new: the Heegaard
$S^2$ in $S^3$, the Heegaard tori in the lens spaces, and the
horizontal splittings.    Each of these has genus different than two
and therefore cannot be isotopic to a Heegaard surface for the knot
exterior, all of which are genus two.

Note that in the filled manifolds, the core of the attached solid
torus is not isotopic into the Heegaard surface $S^2$ for $S^3$, but
is isotopic into both the Heegaard tori for the lens spaces, and the
horizontal splittings, when they occur.

\subsection{N -- The non-level case}

Note that when condition N occurs, the core curve $\gamma$ has some
bridge number $b > 0$ with respect to the Heegaard surface $\Sigma$.
The only filling on a torus knot exterior producing condition N is
the $\frac{1}{0}$--filling which produces $S^3$.  The core $\gamma$
is isotopic into every other Heegaard surface (H or C) for every
other manifold obtained by filling on the torus knot.   Of course,
$\gamma$ is a torus knot in $S^3$ which has bridge number
$b=\min\{p,q\}$ with respect to a sweepout by $S^2$s
(see Schubert \cite{schubert} and Schultens \cite{schultensTorusKnots}).  Tubing $\Sigma$ along the $b$
upper (or lower) bridges, one builds a Heegaard surface for the knot
exterior of genus $g(\Sigma) + b$. Tubing corresponds to
stabilization in the filled manifold, but the resulting surface may
or may not be a stabilized Heegaard surface for the knot manifold.
For  torus knots, this process yields a Heegaard surface of genus
$b$ for the knot exterior, which is irreducible only when its genus
is minimal, that is, when $b=2$.  Such torus knots have $p=2$ or $q=2$
and are excluded for consideration here by the first condition in
\fullref{secOurExamples}.

\subsection{H -- The horizontal case, parallel and boundary stabilization}
\label{secHorizontal}

In case H, we can assume that $\gamma \subset \Sigma$ in the filled
manifold $X(\alpha)$.   The surface $\Sigma^* = \Sigma - N(\gamma)$
is a surface with boundary properly embedded in the knot exterior
$X$.  Moreover, since $\gamma$ is isotopic into the surface, the
meridional slope $\alpha$ intersects the slope $\sigma$ of the
surface $\Sigma^*$ precisely once.  The surface $\Sigma^*$ has two
boundary components and splits $X$ into two handlebodies. Such a
surface is referred to as an {\it almost Heegaard surface} in
Rieck--Sedgwick \cite{rs1}. As noted, $\Sigma$ is a Heegaard surface for every
manifold $X(\alpha')$ where $\alpha' \in L_\sigma$, where $\sigma$
is the slope of $\Sigma^*$. These manifolds differ by Dehn twists in
the curve $\gamma$ in the Heegaard surface $\Sigma$ .

The lens space fillings on the torus knot exterior are a simple example:
$\Sigma$ is the Heegaard torus,  $\Sigma^*$ is the cabling annulus
and its slope is $\sigma = \frac{pq}{1}$.  Moreover, $\Sigma$ is a
Heegaard surface for the entire line $L_{pq/1}$ of lens space fillings.
The collection of Seifert fibered manifolds possessing horizontal
splittings are another example, this time occurring with filling
coefficients $L_{0/1}$.  Note that while most of these splittings are
irreducible, there are a few cases where they reduce ($\frac{1}{0}
\in L_{pq/1}$ and $\frac{1}{0} , \frac{1}{1} \in L_{0/1}$).  It is a
theorem of Casson and Gordon \cite{cgPretzel} that if $\Sigma^*$ is
incompressible then the distance between weakly reducible fillings on
the line is at most 6.  See the appendix of Moriah--Schultens \cite{ms}
for a proof. For horizontal splittings of Seifert fibered spaces, there
are at most two fillings on the line that result in weakly reducible
Heegaard splittings (see Sedgwick \cite{sedgwickSFS}).

\begin{figure}[ht!]
\begin{center}
\labellist\small
\pinlabel {$a$} [l] at 264 418
\endlabellist
\includegraphics[width=.48\hsize]{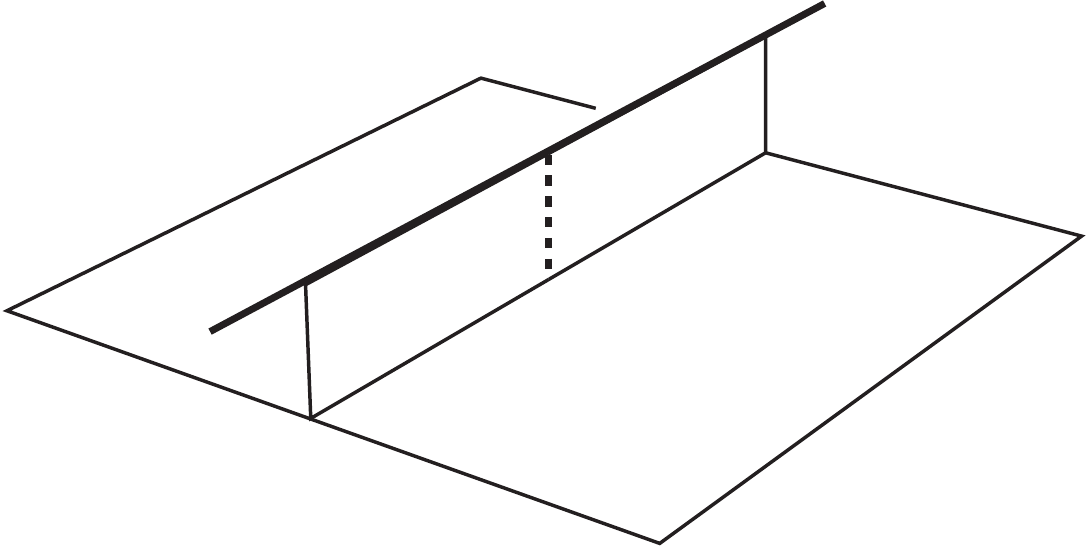} \quad
\includegraphics[width=.48\hsize]{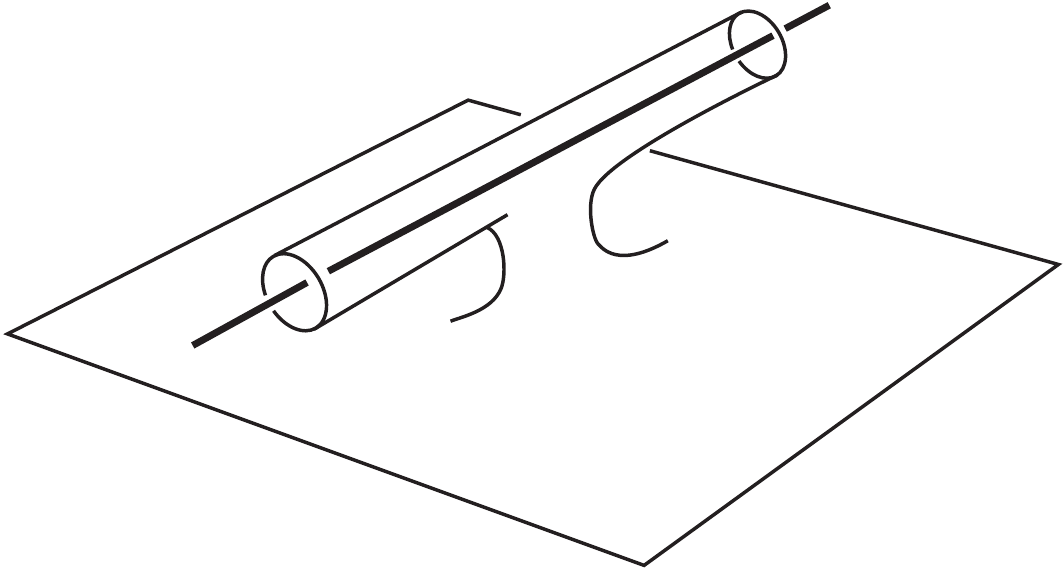}

\caption{Parallel stabilization of a horizonal surface}
\label{figParallelStabilization}
\end{center}
\end{figure}

We can also use $\Sigma$ to form a Heegaard surface for the knot
exterior of genus $g(\Sigma) + 1$ by a process we will call {\it
parallel stabilization}:   Push the surface $\Sigma$ to one side of
the knot, say below, see \fullref{figParallelStabilization}. Now,
$\Sigma$ and $\partial X$ cobound an annulus $A$.  Surgering
$\Sigma$ along $A$ yields $\Sigma^*$, hence $A$ and $\Sigma^*$ have
the same boundary slope on $\partial X$.  Let $a$ be a spanning arc
for $A$.   We can perform {\it parallel stabilization}  on $\Sigma$
by tubing around $\gamma$ and then attaching the tube to $\Sigma$ by
a tube around $a$.   To be more formal, note that the surface
$\partial N(\Sigma \cup a \cup \gamma)$, (neighborhood taken in
$X(\alpha)$) has two components. The parallel stabilization is the
component that has one higher genus than $\Sigma$.  It is a Heegaard
surface for the knot exterior $X$.  Note that we can parallel
stabilize $\Sigma$ in two ways, by starting with  $\Sigma$ either
above or below $\gamma$. These surfaces are not in general isotopic
in $X$, but they are in the filled manifold $X(\alpha)$, because
there they are both stabilizations of the same surface $\Sigma$. It
is not hard to check that parallel stabilization does not depend on
the choice of annulus, even for annuli with different slopes.

\begin{figure}[ht!]
\begin{center}
\labellist\small
\pinlabel {$D$} <-1pt,.5pt> at 248 408
\pinlabel {$D$} <.5pt,.5pt> at 316 440
\endlabellist
\includegraphics[height=1.4in]{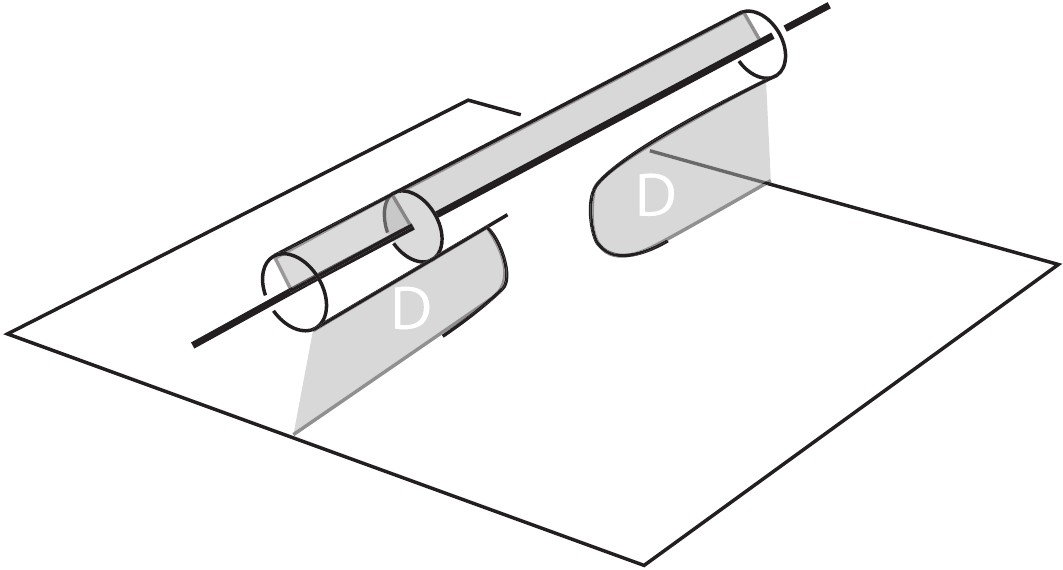}

\caption{The disk $D$ and annuli $A_\alpha$ and $A_\sigma$ in a
parallel stabilization}
\label{figParallelPrim}
\end{center}

\end{figure}

Let $A_\alpha$  be the meridional annulus, $A_\sigma$ be the
longitudinal annulus, and $D$ the disk that appear in \fullref{figParallelPrim}. The annuli $A_\sigma$ and $A_\alpha$ meet in
a single arc, $A_\sigma \cap D = \emptyset$ and $A_\alpha \cap D =
\{pt\}$. Together, $A_\alpha$ and $D$ demonstrate that the parallel
stabilization is $\gamma$--primitive, for $\gamma = \alpha$. But,
$\alpha$ is not unique in this regard. By twisting in $A_\sigma$ we
can construct an $A_{\alpha'}$ meeting $D$ once for any $\alpha' \in
L_\sigma$.   The existence of the triple $(A_\alpha, A_\sigma, D)$
with the specified intersections is  equivalent to the Heegaard
splitting being a parallel stabilization.

\begin{figure}[ht!]
\begin{center}
\includegraphics[width=.48\hsize]{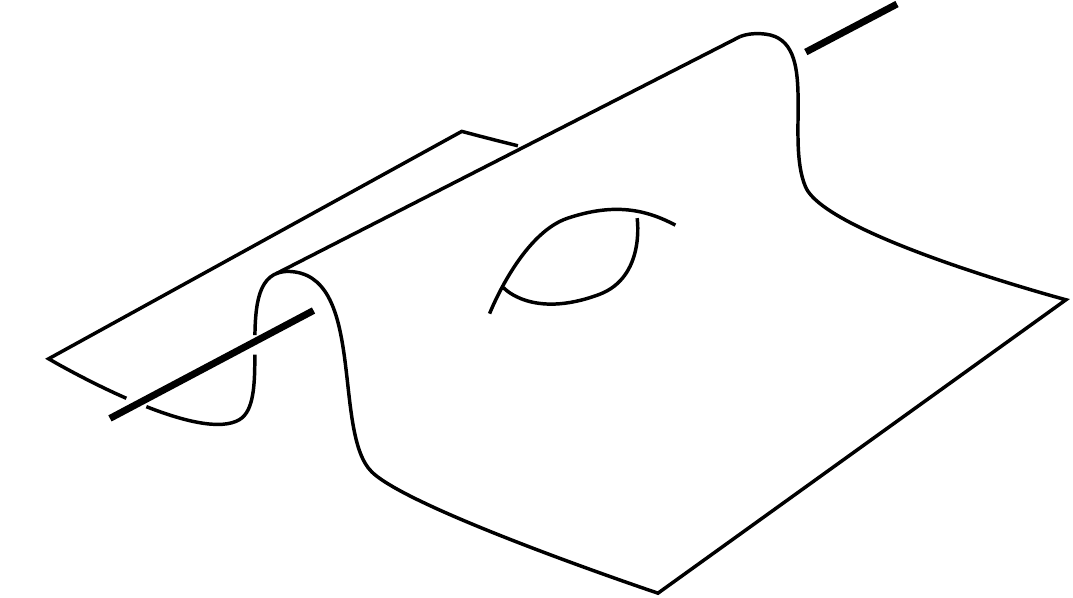}
\includegraphics[width=.48\hsize]{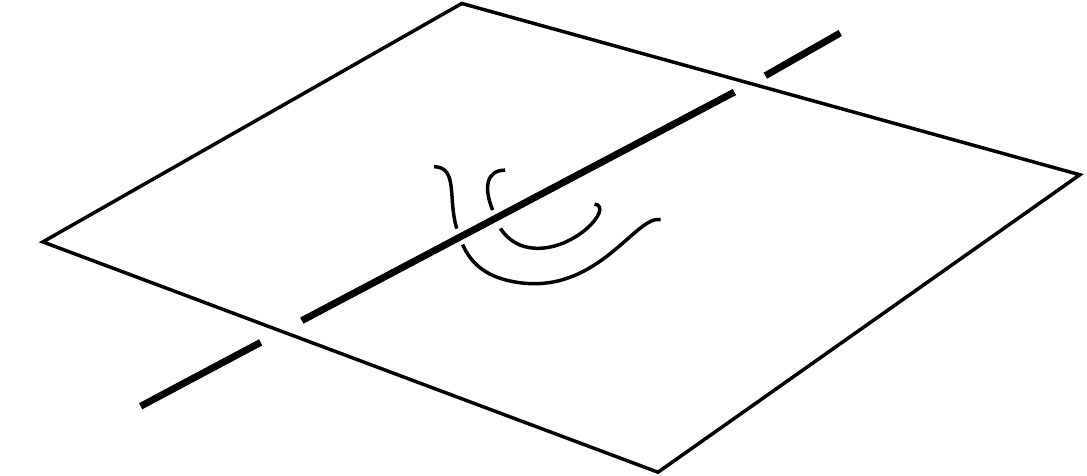}

\caption{Surfaces that are isotopic to the parallel stabilization.
 First make the disk $D$
small, and then flatten the top. This forces the hole to appear as a
tube below the surface.} \label{figParallelIsotopy}
\end{center}
\end{figure}

There is an isotopic picture of the parallel stabilization that will
be very useful in \fullref{secIsotopic}.   See \fullref{figParallelIsotopy}.  Shrink the disk $D$ so that it is small,
and then flatten the top surface. This makes the knot boundary
$\partial X$ appear to be below the surface, and the hole defined by
$D$ now appears to be a tube below the surface. It is also easy to
see $A_\alpha$ as a once punctured disk that meets the disk $D$
once. The annulus $A_\sigma$ runs between the knot and the surface
above it.

If $\Sigma$ was a Heegaard surface for the knot exterior, and not
just a horizontal surface, then we will call the parallel
stabilization of $\Sigma$ a {\it boundary stabilization}. In this
case it can also be viewed as an amalgamation of $\Sigma$ with a
type (ii) splitting of product neighborhood of $\partial X$, see
Moriah \cite{moriahWR}.    Since $\Sigma$ cobounds annuli of every slope
with $\partial X$, and parallel stabilization doesn't depend on this
choice, a boundary stabilization is $\gamma$--primitive for every
slope on $L_\sigma$ for every $\sigma$.   In other words, a boundary
stabilization is $\gamma$--primitive for every slope $\gamma$.

An examination of \fullref{figTunnels} makes it clear that the
inner and outer tunnel systems of the torus knot exterior are
parallel stabilizations of parallel copies of, just above and below,
the standard torus in which $T_{p,q}$ is embedded. However, the
middle splitting is not a parallel stabilization because it  is not
$\gamma$--primitive for any $\gamma$ when it is not isotopic to
either the inner or outer splitting, that is, for $p \not \equiv \pm 1
\pmod q$ and $q \not \equiv \pm 1 \pmod p$ \cite{msTorusKnots}.

When a parallel stabilization is minimal genus then it is clearly
irreducible.  This  is the case with the inner and outer splittings
of the torus knot exterior.  But this need not be the case. In fact,
when filling a torus knot exterior, if the filled manifold has an
irreducible horizontal Heegaard splitting $\Sigma$, then its
parallel stabilization is reducible. Its parallel stabilization
$\Sigma'$ is a splitting of the knot manifold, a Seifert fibered
space with boundary. As already noted it is proved by Schultens
\cite{schultensSFSBoundary} that such splittings are vertical and
irreducible only when they are minimal genus.  But if $p > 3$ then
the genus of $\Sigma'$ is $2(p-1)(q-1)+1 > 3$, strictly greater than
two, the genus of the knot exterior.

\section{New Heegaard surfaces -- the results}
\label{secResults}

We now survey results that restrict  (N) and (H).

\subsection{Moriah and Rubinstein}

In order to prove the existence of knots with super-additive tunnel
number, Moriah and Rubinstein \cite{mr} needed to show that there
are fillings on a knot exterior for which the Heegaard genus does
not drop at all. If the genus does drop, then the filled manifold
has a Heegaard surface of genus lower than that of every Heegaard
surface for the knot exterior, that is, the filled manifold possesses a
new Heegaard surface.

In fact, they conclude that there {\it is} a finite list of
candidates of bounded genus for all but a finite number of fillings,
and the genus does not degenerate for ``most'' fillings, as is
demonstrated by the following theorem and corollary. Their theorem
applies to manifolds with multiple torus boundary components.  Here
we have restated it only for knot--manifolds and in terms consistent
with our discussion:

\begin{thm}[Moriah and Rubinstein -- rephrased]

Let $X$ be hyperbolic knot--manifold and $g$ a positive integer.
Then there is finite set of slopes $\mathcal N_X$ and a collection
of Heegaard and horizontal surfaces $\Sigma_1, \ldots ,\Sigma_k
\subset X$, so that if $\alpha \not \in \mathcal N_X$ then any
Heegaard splitting of $X(\alpha)$ of genus less than or equal to $g$
is isotopic to one of the $\Sigma_i$.
\end{thm}

Their theorem implies several things for Heegaard surfaces of
bounded genus, in particular minimal genus, in hyperbolic
knot--manifolds. Condition N occurs for at most finitely many
slopes, those in $\mathcal N_X$. Away from this finite number of
slopes, the core is isotopic into {\it every} Heegaard surface of
bounded genus. And, H occurs for all slopes on a finite (possibly
empty) set of lines $L_{\beta_1}, \ldots, L_{\beta_k}, k \geq 0$.

By avoiding a finite set of slopes for N and a finite set of lines
for H, {\it every} Heegaard surface of bounded genus for the filled
manifold is isotopic to a Heegaard surface for the knot exterior. In
particular, the Heegaard genus of these manifolds is the same as
that of the knot exterior.  Of course, if the genus decreases and we
are not in situation N, then the discussion on the horizontal case H
(\fullref{secHorizontal}) shows that the genus decreases by at
most one.

A corollary of their theorem is:

\begin{cor}[Moriah and Rubinstein]

Let $X$ be a hyperbolic knot--manifold.  Then there exists a finite
set of slopes $\mathcal N_X$ and a finite set of lines $\mathcal
H_X$ so that:

\begin{enumerate}
\item $\alpha \notin \mathcal N_X \cup \mathcal H_X \implies g(X(\alpha)) = g(X)$
\item $\alpha \notin \mathcal N_X \implies g(X(\alpha)) \geq g(X)-1$.
\end{enumerate}
where $g$ denotes the Heegaard genus of the manifold.
\end{cor}

This theorem forbids new Heegaard splittings of bounded genus for
most filled manifolds, but leaves open that possibility for new
Heegaard surfaces when we do not bound genus. That is, as we
increase the bound $g$ the sets of slopes for conditions N and H
could grow in an unbounded fashion.   In fact, as discussed in
\fullref{rs}, this does not happen.

\subsection{Rieck}

Rieck \cite{rieckDehnFilling} took a topological approach to the
same problem and computed numeric bounds on the distance between
``bad'', that is, type N, fillings. Suppose that N occurs with respect
to Heegaard surfaces $\Sigma_1 \subset X(\alpha_1)$ and $\Sigma_2
\subset X(\alpha_2)$. Since the core is not isotopic into either
$\Sigma_1 \subset X(\alpha_1)$ or $\Sigma_2 \subset X(\alpha_2)$, it
can be put into non-trivial thin position with respect to sweepouts
by each surface. Furthermore, thick level surfaces for each
sweepout, when regarded as punctured surfaces in the knot exterior,
will intersect essentially. This approach yields a bound on the
distance between the slopes $\alpha_1$ and $\alpha_2$:

\begin{thm}[Rieck -- rephrased]
Let $X$ be an anannular knot manifold.   Suppose that the core of
the attached solid torus is not isotopic into Heegaard surfaces
$\Sigma_1 \subset X(\alpha_1)$ and $\Sigma_2 \subset X(\alpha_2)$.
Then $\Delta(\alpha_1,\alpha_2) < 18g_1g_2 + 18g_1 + 18g_2 + 18$,
where $g_1$ and $g_2$ are the genera of $\Sigma_1$ and $\Sigma_2$,
respectively.
\end{thm}

In \cite{rieckGenusReducing}, Rieck examines the relationship
between Dehn filling and Heegaard structure from a
 different viewpoint.  He asks which manifolds possess a {\it genus
 reducing knot}, a knot for which infinitely many surgeries decrease
 Heegaard genus.     He answers this question for all totally
 orientable Seifert fibered spaces other than those with base space
 $S^2$ and three or fewer exceptional fibers.  Almost all of the
 considered Seifert fibered spaces do contain contain genus reducing
 knots, the exception being Seifert fibered spaces possessing a
 horizontal Heegaard surface of one of two special types.

\subsection{Rieck and Sedgwick}
\label{rs}

Condition H is explored further in Rieck and Segwick \cite{rs1}.  As we have already
noted, in that case we can form the surface $\Sigma^* = \Sigma -
N(\gamma)$, a properly embedded surface in the knot--manifold $X$.

\begin{thm}[Rieck and Sedgwick]
Suppose that the core of the attached solid torus is isotopic into
$\Sigma$ a Heegaard surface for a filled manifold $X(\alpha)$. Then
one of the following holds:

\begin{enumerate}
\item $\Sigma$ is a Heegaard surface for $X$ (perhaps after an
isotopy in $X(\alpha)$, or,
\item the slope of the almost Heegaard surface $\Sigma^*$ is the
boundary slope of a separating essential surface of genus less than
or equal to that of $\Sigma^*$.
\end{enumerate}
\end{thm}

If the second conclusion occurs, then the slope $\alpha$ is one that
intersects the slope of an essential surface exactly once.   Hatcher
has shown such slopes to be finite in number \cite{hatcher}, so if H
occurs we know that the slope $\alpha$ belongs to one of a finite
number of lines defined by slopes of essential surfaces. This
improves the earlier work of Moriah and Rubinstein because the
knot--manifold is not required to be hyperbolic, it applies to all
surfaces without a bound on genus, and a connection is made between
these slopes and the slopes of essential surfaces.

It would be nice if the surface $\Sigma^*$ were itself an essential
surface. However, the method of proof is similar to that of Casson
and Gordon \cite{cg} and may require that $\Sigma^*$ is modified by
compressions and annulus swaps to obtain an essential surface.

In their second paper \cite{rs2} Rieck and Sedgwick continued their
investigation into the Heegaard structure of filled manifolds. We
can assume that the knot--manifold $X$ is given via a one-vertex
triangulation, see Jaco and Sedgwick \cite[Theorem~3.2]{jacoSedgwick}.
They then prove (slightly rephrased):

\begin{thm}[Rieck and Sedgwick]
Let $\mathcal T$ be a one-vertex triangulation of the knot--manifold
$X$.   If $\gamma$, the core of the attached solid torus, is not
isotopic in $X(\alpha)$ into a Heegaard surface $\Sigma$, then the
slope $\alpha$ is either the slope of a boundary edge of the
triangulation or the slope of a normal or almost normal slope in
$(X,\mathcal T)$.
\end{thm}

The proof of the theorem follows Thompson's proof \cite{thompsonS3}
that a triangulation of $S^3$ contains an almost normal $S^2$. The
1--skeleton of the triangulation is put in thin position with
respect to a sweepout given by the Heegaard surface $\Sigma$.    If
a boundary edge is isotopic into the Heegaard surface, we either
have an edge with slope $\alpha$ or the core $\gamma$ is isotopic
into $\Sigma$.   Otherwise, we are able to find a non-trivial thick
level,  yielding a normal or almost normal surface with slope
$\alpha$ in $\Sigma$.   Bachman \cite{bachmanANBoundary} has a
similar result.

They then apply a theorem of Jaco and Sedgwick \cite{jacoSedgwick}
stating that there are only finitely many slopes bounding normal and
almost normal surfaces in such a triangulation.

\begin{thm}[Jaco and Sedgwick]
Let $X$ be a knot-manifold with a triangulation $\mathcal T$ that
restricts to a one--vertex triangulation on $\partial X$.   Then
there are only a finite number of slopes realized as the slopes of
embedded normal and almost normal surfaces in $(X,\mathcal T)$.
\end{thm}

The proof is an analog of Hatcher's proof \cite{hatcher} that
there are a finite number of slopes bounding essential surfaces in a
knot--manifold.  It is shown that normal or almost normal surfaces
that are compatible, meaning that their normal sum is well defined,
must have the same slope or their sum produces trivial curves in the
boundary. Whereas Hatcher's theorem relies on Floyd and Oertel's
work with branched surfaces \cite{floydOertel}, this proof appeals
to similar properties of normal and almost normal surfaces in a
triangulation.

\begin{cor}[Rieck and Sedgwick, rephrased]
Let $X$ be a knot--manifold.  Then there exists a finite set of
slopes $\mathcal N_X$ in $\partial X$ so that if $\alpha \notin
\mathcal N_X$, then the core of the attached solid torus $\gamma$ is
isotopic into every Heegaard surface for $X(\alpha)$.
\end{cor}

Note that this theorem does not require a bound on genus and applies
to non-hyperbolic as well as hyperbolic knot exteriors.

\subsection{Summary of known results}

For clarity, we offer a summary of the known results.   Recall the
trichotomy offered at the start of this section.  The core of the
attached solid torus $\gamma$ is either (N)ot Level, (H)orizontal
but not a core, or a (C)ore of a given Heegaard surface $\Sigma$.

\begin{thm}
\label{thmSummary}
 Let $X$ be a knot--manifold.  Then there is a
finite set of slopes $\mathcal N_X$ and a finite set of lines
$\mathcal H_X$ so that:

\begin{enumerate}
\item If $\alpha \notin \mathcal N_X$, then the core of the attached
solid torus is isotopic into every Heegaard surface for $X(\alpha)$,
and in particular $g(X) - 1 \leq g(X(\alpha)) \leq g(X)$, and,
\item If $\alpha \notin \mathcal N_X \cup \mathcal H_X$ then $X(\alpha)$
does not contain a new Heegaard surface, that is, every Heegaard
surface for $X(\alpha)$ is isotopic (in $X(\alpha)$) to a Heegaard
for $X$, and in particular $g(X(\alpha)) = g(X)$.
\end{enumerate}
\end{thm}

\begin{rmk}
For fillings on any torus knot exterior we can take $\mathcal N_X =
\{\frac{1}{0}\}$ and $\mathcal H_X = \{L_{pq/1} \cup L_{0/1}\}$.
\end{rmk}

\section{What is not known?}

While we have answered many of the questions regarding the Heegaard
structure of filled manifolds, there are at least several that
remain.

\subsection{Destabilization}
\label{secDestab}

\begin{quest}
Let $\Sigma$ be a non-stabilized  Heegaard surface for $X$.  What
can be said about the set of the fillings for which $\Sigma$
destabilizes?
\end{quest}

When filling a torus knot exterior the  inner, outer and middle
splittings necessarily destabilize when the obtained manifold is a
lens space or $S^3$, that is, for the fillings on the line $L_{pq/1}$.

This is evident in the fact that each of these splitting is padded.
A Heegaard surface $\Sigma$ is said to be {\it padded} if there
exists a triple $(P,A_\sigma,D)$ where:

\begin{enumerate}
\item  $P$ is a punctured disk in the compression body, (a {\it
punctured disk} is a planar surface in a compression body that has
one boundary component, its ``boundary'', in the Heegaard surface
and all others, the ``punctures'', in $\partial X$)

\item $A_\sigma$ is a vertical annulus in the compression body with
slope $\sigma$ on $\partial X$.

\item $D$ is a disk in the handlebody,

\item $|\partial P \cap
\partial D| =1$

\item $|\partial A_\sigma \cap
\partial D |= 0$.
\end{enumerate}

Condition (5) clearly implies that padded splittings are weakly
$\sigma$--primitive. But in fact, it is much stronger:  For a
strongly irreducible Heegaard surface the results of
Culler--Gordon--Luecke--Shalen \cite{cgls} and
Wu \cite{wu} imply that $\sigma$ meets the slope of the punctures once.
Twisting $P$ in the annulus $A_\sigma$ yields a destabilizing pair
$(P_\alpha,D)$ for every slope $\alpha \in L_\sigma$.

We have observed that the inner and outer splittings $\Sigma_i$ and
$\Sigma_o$ are parallel stabilizations.  This is a very strong form
of padded where the punctured disk $P$ is actually an annulus (take
$P = A_\alpha, A_\sigma = A_\sigma, D = D$).  The middle splitting
$\Sigma_m$ is also padded where $\sigma = \frac{pq}{1}$, but it is
slightly harder to see. The handlebody component of $X - \Sigma_m$
is a regular neighborhood of the union of the critical fibers joined
by a spanning arc $a_{io}$ for the annulus $A_{io}$ running between
them, that is, $N(f_i \cup a_{io} \cup f_o)$.  The arc $a_{io}$
punctures the cabling annulus $A$ once so we may think of the
compression body as a product $T' \times [-1,1] - N(T_{p,q})$, where
$T'$ is a copy of $T$ with a single puncture. Take $D$ to be the
meridional disk for the inner solid torus, $A \subset V$ as the
vertical annulus $T_{p,q} \times [0,1]$ with slope $\frac{pq}{1}$
and boundary towards the outside solid torus, and $P = \{a\} \times
I - N(T_{p,q})$, where $a$ is any arc properly embedded in $T'$ that
meets $\partial D$ exactly once.

It may be surprising to realize that there is a non-stabilized
Heegaard splitting that destabilizes for every filling. As noted in
\fullref{secNewSplittings}, a boundary stabilization is
$\gamma$--primitive for every $\gamma$ so it necessarily
destabilizes for every filling. Hence the question becomes: Can a
boundary stabilized Heegaard splitting be non-stabilized in the
knot-manifold? Yes, as is shown in Sedgwick \cite{sedgwickGenus2} and
Moriah--Segwick \cite{msBS}. Let
$\Sigma_{ij}$ be one of the genus two splittings for $\{\textit{pair
of pants}\} \times S^1 $ discussed in \fullref{secPOP}.  Then
the annuli $A_{ik}$ and $A_{jk}$ are annuli that we can use to
parallel stabilize (boundary stabilize, since we start with a
Heegaard surface).   In each case, we obtain the genus three
splitting $\Sigma_{123}$ which is identified by the partition of
boundary components $\{\partial X_1,
\partial X_2, \partial X_3
\parallel \emptyset\}$.  It is induced by the pair of arcs $a_{12} \cup a_{23}$
(or $a_{12} \cup a_{13}$ or $a_{13} \cup a_{23}$).    There are also
examples of boundary stabilized but non-stabilized Heegaard
splittings for manifolds with two boundary components \cite{msBS}.
The question remains open for knot-manifolds.

In general, we expect to be able to say little about a Heegaard
splitting that destabilizes for a finite number of fillings. But
what if it destabilizes for infinite number of fillings?  The above
examples demonstrate that padded splittings destabilize for
infinitely many fillings and a boundary stabilization destabilizes
for all fillings. We conjecture the converses:

\begin{cnj} Suppose that a non-stabilized Heegaard surface
destabilizes for infinitely many fillings.  Then it is padded.
\end{cnj}

\begin{cnj}Suppose that an irreducible Heegaard surface
destabilizes for all fillings.  Then it is a boundary stabilization.
\end{cnj}

It might appear that \fullref{thmSummary} would answer the
questions raised in this section. But \fullref{thmSummary} only
restricts the set of slopes for which a given Heegaard surface
destabilizes to a {\it new} Heegaard surface.   But, it is possible
for a non-stabilized Heegaard surface to destabilize in the filled
manifold where it becomes isotopic to another Heegaard surface for
the knot exterior.  For example, a boundary stabilization of a
Heegaard surface $\Sigma$ destabilizes in every filling, and the
destabilized surface is isotopic to $\Sigma$.

But, more importantly, \fullref{thmSummary} doesn't say anything
about the structure of the surface destabilizing, even if it is of
minimal genus:   Let $\Sigma$ be a minimal genus Heegaard surface
for $X$ that destabilizes in infinitely many fillings, in particular
in some filled manifold $X(\alpha)$ where $\alpha \notin \mathcal
N_X$. Then $\Sigma$ destabilizes to a Heegaard surface $\Sigma'$ for
$X(\alpha)$ and the core of the attached solid torus is isotopic
into $\Sigma'$. We can construct $\Sigma''$, a parallel
stabilization of $\Sigma'$, that is a Heegaard surface for $X$ of
the same genus as that of $\Sigma$.  The parallel stabilization
$\Sigma''$ destabilizes in $X(\alpha)$, as well as in an infinite
number of other filled manifolds.  But we don't know that $\Sigma$
and $\Sigma''$ are the same surface in $X$, all we know is that they
become isotopic in $X(\alpha)$.

In fact, the above situation is not vacuous.   The middle splitting of the
$(p,q)$--torus knot exterior is padded but not a parallel stabilization.
(A parallel stabilization is $\gamma$--primitive for all slopes meeting
the slope of the annulus $A$ once, but the middle splitting is not
$\gamma$--primitive for any $\gamma$ as proved in Moriah--Sedgwick
\cite{msTorusKnots}.)     Not only does it destabilize for fillings on
the line $L_{pq/1}$, it also becomes isotopic to the inner and outer
splittings for fillings on the same line.

\subsection{Isotopic Surfaces}
\label{secIsotopic}
 The above discussion leads us to our next question:

 \begin{quest}
 Suppose that $\Sigma_1$ and $\Sigma_2$ are non-isotopic Heegaard surfaces for $X$. In
which fillings do $\Sigma_1$ and $\Sigma_2$ become isotopic?
\end{quest}

This question is also of interest if we wish to understand what
happens to the oriented Heegaard tree $\mathcal H_X^\pm$ after
filling. We know that some vertices are identified in every filled
manifold.

\begin{figure}[ht!]
\begin{center}
\includegraphics[height=0.5in]{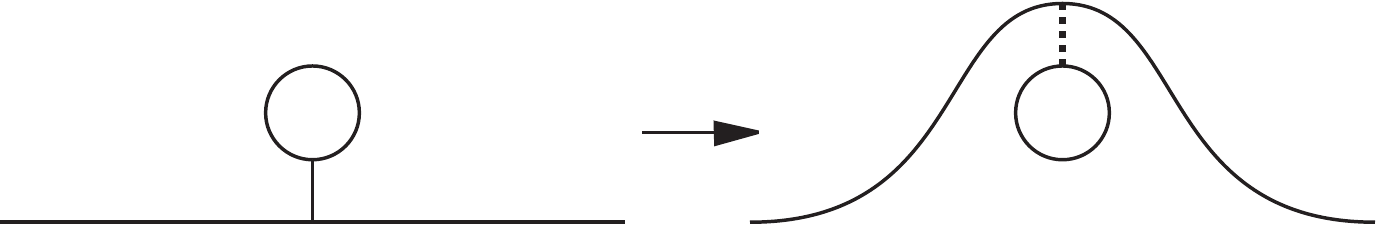}
\caption{An annulus swap on a horizontal surface}
\label{figAnnulusSwap}
\end{center}
\end{figure}

An annulus swap is a procedure that allows us to find, in the
knot-manifold, surfaces that become isotopic after filling.  See
\fullref{figAnnulusSwap}. Suppose that $\Sigma$ is a horizontal
surface in $X$,  so that there is an annulus $A_\sigma$ with one
boundary component in $\Sigma$ and another a non-trivial curve of
slope $\sigma$ in $\partial X$.  Then $\partial N(\Sigma \cup
A_\sigma \cup
\partial X)$, where the neighborhood is taken in $X$, consists of
two surfaces, one isotopic in the neighborhood to $\Sigma$ and the
other not, call it $\Sigma'$. Then we say that $\Sigma'$ is obtained
from $\Sigma$ by an {\it annulus swap} along $A_\sigma$.   There is
an annulus $A_\sigma'$ running from $\Sigma'$ to $\partial X$, also
with slope $\sigma$, that we can use to reverse the operation,
obtaining $\Sigma$ from $\Sigma'$ by swapping $\Sigma'$ along
$A_\sigma'$. More importantly, after filling along any slope $\alpha
\in L_\sigma$ the surfaces $\Sigma$ and $\Sigma'$ co--bound a
product, $\Sigma \times I$, and are therefore isotopic in
$X(\alpha)$.   We will say that $\Sigma$ and $\Sigma'$ are {\it
$\sigma$--swap equivalent} if $\Sigma'$ can be obtained from
$\Sigma$ by a swap along an annulus $A_\sigma$ with slope $\sigma$.
They are {\it swap equivalent} if they are $\sigma$--swap equivalent
for some $\sigma$.   And, we will say that $\Sigma$ and $\Sigma'$
are {\it weakly swap equivalent} if there is a sequence of swaps
(with no restriction on slopes) taking $\Sigma$ to $\Sigma'$.

Distinct Heegaard surfaces can be swap equivalent.  For example, the
dual tunnels of Morimoto and Sakuma \cite{morimotoSakuma} are
$\frac{1}{0}$--swap equivalent.  The following lemma shows that this
only happens under specific circumstances:

\begin{lem}
Suppose that $\Sigma$ is a Heegaard surface for $X$.  Then $\Sigma$
is $\gamma$--swap equivalent to a Heegaard surface $\Sigma'$ for $X$
if and only if $\Sigma$ is $\gamma$--primitive.
\end{lem}

\begin{proof}

Let $V \cup_\Sigma W$ and $V' \cup_{\Sigma'} W'$ be the
decompositions induced by $\Sigma$ and $\Sigma'$, where $W$ and $W'$
are handlebodies and $V$ and $V'$ are the components containing
$\partial X$. We know that $V$ is a compression body and must decide
when $V'$ is as well.

Let $N$ be a regular neighborhood, $N = N(A_\gamma \cup \partial
X)$, taken in $V$.   Its boundary $\partial N$ consists of two
2--tori $T_1 \cup T_2$, one of which, say $T_2$, is $\partial X$.
The torus $T_1$ is composed of two annuli $T_1 = A_V \cup A_\Sigma$
where $A_\Sigma \subset \Sigma$ is a regular neighborhood in
$\Sigma$ of the curve component of $\partial A_\gamma \setminus
\partial X$ and $A_V = T_1 \setminus A_\Sigma$ is properly embedded in $V$.

We obtain $V'$ by gluing $N$ to $W$ along $A_\Sigma$, that is, $V' = W
\cup_{A_\Sigma} N$.  Since, $N$ is a homeomorphic to a product,
$\textit{2--torus} \times I$, the component $V'$ is a compression
body if and only if there is an essential disk $D \subset W$ that
meets $A_\Sigma$ in a single essential arc.  This occurs if and only
if $D$ meets the boundary of the annulus $A_\gamma$ in a single
point, that is, $\Sigma$ is $\gamma$--primitive.
\end{proof}

It follows that if the Heegaard surfaces $\Sigma$ and $\Sigma'$ are
swap equivalent but one of them is not $\gamma$--primitive for some
$\gamma$, then the sequence of swaps must have length at least two.
This still places restrictions on $\Sigma$ and $\Sigma'$ as detailed
in the following lemma.  Note that $\Sigma$ is not assumed to be a
Heegaard surface, although if it is swap equivalent to a Heegaard
surface, it must be a horizonal surface.

\begin{lem}
\label{lemSwapWGammaPrim} Suppose that we can perform two distinct
annulus swaps on $\Sigma$. Swapping along $A_{\sigma_1}$ yields
$\Sigma_1$, and swapping along $A_{\sigma_2}$ produces $\Sigma_2$,
where $A_{\sigma_1}$ and $A_{\sigma_2}$ are not isotopic.   If
either $\Sigma_1$ or  $\Sigma_2$ is a Heegaard surface, then
$\Sigma_1$ is weakly $\sigma_1$--primitive and $\Sigma_2$ is weakly
$\sigma_2$--primitive.
\end{lem}

\begin{proof}
Isotope $A_{\sigma_1}$ and $A_{\sigma_2}$  to intersect minimally.
It follows that neither has inessential arcs of intersection with
endpoints in $\partial X$.   We claim that we can find essential
disks $D_1$ and $D_2$ (possibly the same) that compress $\Sigma$
towards $\partial X$ and whose boundaries are disjoint from
$A_{\sigma_1}$ and $A_{\sigma_2}$, respectively. This would prove
the lemma.

If $A_{\sigma_1} \cap A_{\sigma_2}$  contains an inessential arc,
then it is inessential in both annuli, and its endpoints are on
$\Sigma$. An outermost inessential arc on $A_{\sigma_2}$ bounds a
disk in $A_{\sigma_2}$ and a disk in $A_{\sigma_1}$.  The union of
these disks is a disk $D_1$ that is essential (minimality) and
disjoint from $A_{\sigma_1}$.   A symmetric argument yields a disk
$D_2$ disjoint from $A_{\sigma_2}$. We can therefore assume that the
annuli meet only in essential arcs, so they are cut into rectangles.

If $\sigma_1$ and $\sigma_2$ do not have the same slope, then we can
form a collection of disks that are shaped like a box without a top,
each disjoint from both annuli.  The four sides of the box are the
sub-rectangles of the annuli and the bottom of a box is a rectangle
in $\partial X$. Furthermore, it is not possible for all of our
boxes to be inessential disks, this would imply that $\Sigma$ is
peripheral in $X$ and, in turn, show that $X$ is a solid torus, not
a proper knot--manifold.

If $\sigma_1$ and $\sigma_2$ do have the same slope then the annuli
are disjoint. In that case, we work in some $X(\alpha), \alpha \in
L_{\sigma_1}$.   Then the annulus $A = A_{\sigma_1} \cup
A_{\sigma_2} \subset X(\alpha)$ is a properly embedded annulus in a
handlebody ($\Sigma_1$ is a Heegaard surface).   It can't be
peripheral as this would imply that $A_{\sigma_1}$ and
$A_{\sigma_2}$ are isotopic in $X$.    So, it boundary compresses to
a disk $D$ that is essential and disjoint from both $A_{\sigma_1}$
and $A_{\sigma_2}$.
\end{proof}

Suppose that $\Sigma$ and $\Sigma'$ are equivalent after a sequence
of swaps of length $n$, but not for a sequence of shorter length:
$\Sigma = \Sigma_0 \leftrightarrow \Sigma_1 \leftrightarrow \cdots
\leftrightarrow \Sigma_n = \Sigma'$. If $n=1$ $\Sigma$ and $\Sigma'$
are both $\gamma$--primitive, hence weakly $\gamma$--primitive, for
some $\gamma$. If $n >1$, then the previous lemma can be applied
with $\Sigma = \Sigma_1$ (and $\Sigma = \Sigma_{n-1}$) to show that
they are both weakly $\gamma$--primitive.   It follows that only
weakly $\gamma$--primitive Heegaard surfaces are swap equivalent to
other Heegaard surfaces.

\begin{cor}
Suppose that $\Sigma$ and $\Sigma'$ are swap equivalent Heegaard
surfaces.  Then each is weakly $\gamma$--primitive for some
$\gamma$.
\end{cor}

\begin{figure}[ht!]
\begin{center}
\includegraphics[width=.48\hsize]{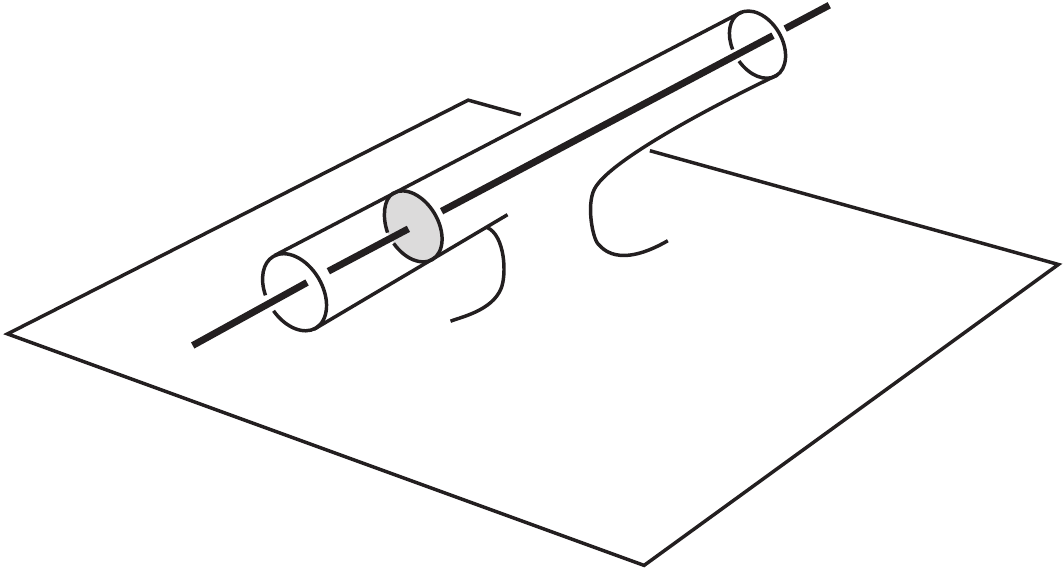}
\includegraphics[width=.48\hsize]{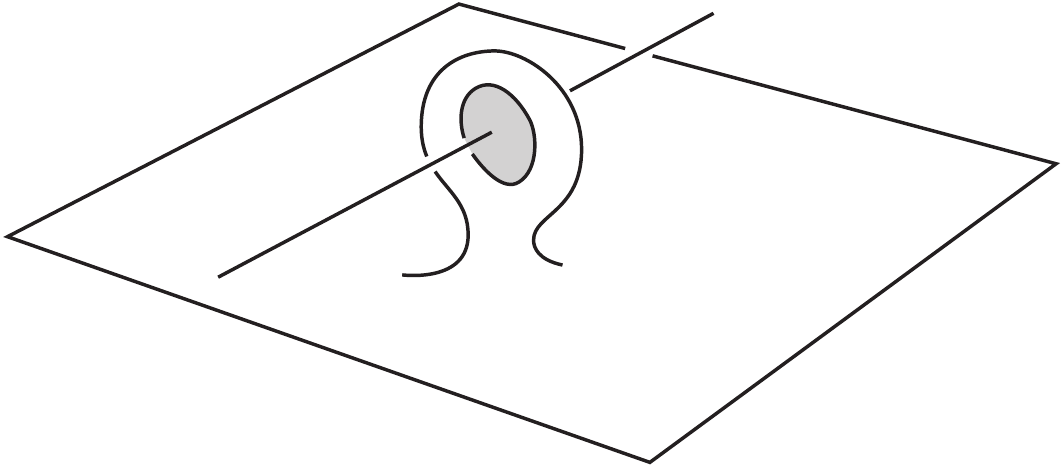}

\caption{Parallel stabilizations that are swap equivalent along
meridional annuli}
\label{figParallelSwap}
\end{center}
\end{figure}

It is particularly interesting to see what happens when we perform a
swap along a parallel stabilization. This case is covered by the
following lemma:

\begin{lem}
Let $\Sigma$ and $\Sigma'$ be horizontal surfaces that are related
by a swap along an annulus with slope $\sigma$. Then the parallel
stabilizations of $\Sigma$ and $\Sigma'$ are swap-equivalent for
every slope $\alpha \in L_\sigma$.
\end{lem}

\begin{proof}
Choose a swap annulus $A_\alpha$ for the parallel stabilization of
$\Sigma$ where $\alpha \in L_\sigma$. \fullref{figParallelSwap}
indicates the situation in $X(\alpha)$, where $A_\alpha$ appears to
be meridional.   The parallel stabilization on the left and the
surface on the right are swap equivalent along the meridional annuli
pictured. It is easier to see the swap from right to left. In that
case the swap just tubes the surface towards and then along the
knot, yielding the surface and dual annulus on the left.   But the
surface on the right hand side is isotopic to an upside--down copy
of the surface in \fullref{figParallelIsotopy}. Thus, the right
hand surface is a parallel stabilization of $\Sigma'$, the surface
obtained if we had swapped $\Sigma$ before parallel stabilizing.
\end{proof}

Recall that if two surfaces are swap equivalent for some slope
$\sigma$, then the surfaces are isotopic in $X(\alpha)$ for every
$\alpha \in L_\sigma$.

\begin{cor}
\label{corIsotopicLL} Let $\Sigma$ and $\Sigma'$ be horizontal
surfaces that are related by a  swap along an annulus with slope
$\sigma$.  Then the parallel stabilizations of $\Sigma$ and
$\Sigma'$ are isotopic in $X(\alpha)$ for every $\alpha \in
LL_\sigma$.
\end{cor}

\begin{rmk}
\label{rmkSign} Since the isotopy passes through the (filled) knot,
it changes the sign of an oriented surface (see \fullref{secStabTree}).
\end{rmk}

A parallel stabilization can also be swapped along the annulus with
slope $\sigma$.  The proof of the following lemma is just the
observation that a parallel stabilization, see \fullref{figParallelStabilization}, is just a regular stabilization
followed by an annulus swap along the slope $\sigma$.

\begin{lem}
\label{lemSigmaSwap} Suppose that $\Sigma$ is a horizontal surface
with annulus $A_\sigma$ running from $\Sigma$ to $\partial X$. Then
the parallel stabilization of $\Sigma$ is isotopic to a
stabilization of $\Sigma$ in $X(\alpha)$ for any $\alpha \in
L_\sigma$.
\end{lem}

\begin{rmk}
If the horizonal surface is not a Heegaard surface, then its
stabilization is also a horizontal surface that is not a Heegaard
surface. (Stabilization and destabilization do not convert a
Heegaard to a non-Heegaard surface.)

\end{rmk}

We now return our focus to the fillings on the $(p,q)$--torus knot
exterior and show that annulus swaps explain when any pair of genus
two Heegaard surfaces for the knot exterior become isotopic.  (We
will continue to restrict the class of torus knots as in \fullref{secOurExamples}).

The inner and outer splittings become isotopic, but remain distinct
from the middle splitting, in fillings on the line of lines
$LL_{pq/1}$. This is explained by the lemma and corollary above.
First note that the Heegaard tori $T_i =
\partial N(f_i)$ and $T_o = \partial N(f_o)$ are horizontal surfaces
that are swap equivalent via an annulus with slope $\frac{pq}{1}$,
hence they are isotopic in the lens space fillings on the line
$L_{\unfrac{pq}{1}}$. The isotopy takes $T_i^\pm$ to $T_o^\mp$ and
vice-versa. Since the Heegaard surfaces $\Sigma_i$ and $\Sigma_o$
are their respective parallel stabilizations, they are isotopic for
all slopes $\beta \in LL_{\unfrac{pq}{1}}$ by \fullref{corIsotopicLL}.   As oriented surfaces, this isotopy takes
$\Sigma_i^\pm$ to $\Sigma_o^\mp$, and changes the partition in $X$
from $\{f_i, \partial X
\parallel f_o\}$ to $\{f_i \parallel \partial X,f_o\}$.

\begin{figure}[ht!]
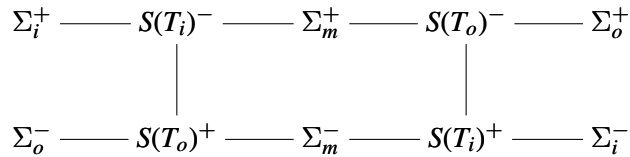

$$\bfig
  \barrsquare/-``-`-/<500,400>[\Sigma_i^+`S(T_i)^-`\Sigma_o^-`S(T_o)^+;```]
  \barrsquare(500,0)/-```-/<500,400>[\phantom{S(T_i)^-}`\Sigma_m^+`
    \phantom{S(T_o)^+}`\Sigma_m^-;```]
  \barrsquare(1000,0)/-```-/<500,400>[\phantom{\Sigma_m^+}`S(T_o)^-`
    \phantom{\Sigma_m^-}`S(T_i)^+;```]
  \barrsquare(1500,0)/-`-``-/<500,400>[\phantom{S(T_o)^-}`\Sigma_o^+`
    \phantom{S(T_i)^+}`\Sigma_i^-;```]
  \efig$$
\caption{Genus two surfaces that are equivalent via swaps along
annuli with slope $\frac{pq}{1}$.  These surfaces are all isotopic
in the lens spaces $X(\alpha)$ where $\alpha \in L_{pq/1}$.}
\label{figPQSwaps}

\end{figure}

All three genus two Heegaard surfaces become isotopic along the line
$L_{\unfrac{pq}{1}}$ of lens space fillings. \fullref{lemSigmaSwap}
shows that the inner and outer surfaces $\Sigma_i^+$ and
$\Sigma_o^+$ become isotopic to stabilizations of the horizontal
inner and outer Heegaard tori $S(T_i^-)$ and $S(T_o^-)$,
respectively.  Note that while $S(T_i^-)$ and $S(T_o^-)$ are not
Heegaard surfaces, they are $\frac{pq}{1}$--primitive.  The
handlebody for the middle splitting $\Sigma_m$ is a regular
neighborhood of the union of the exceptional fibers $f_i$ and $f_o$
joined by a spanning arc for the annulus $A_{io}$ running between
them.  So we can find two annuli from $\partial X$ to $\Sigma_m$,
each with slope $\frac{pq}{1}$.  One runs between $\partial X$ and
$T_i$, the other between $\partial X$ and $T_o$. Swapping along
these annuli changes $\Sigma_m^+$ to the stabilized horizontal
surfaces $S(T_i^-)$ and $S(T_o^-)$, respectively.  These
observations yield the swap diagram in \fullref{figPQSwaps}, all
using annuli with slope $\frac{pq}{1}$. In particular, it explains
how all three surfaces flip in the lens spaces.   We suspect that
this swap diagram is complete.

Annulus swaps explain every post--filling isotopy of Heegaard
surfaces for the torus knot exteriors we consider.  This leads to
several questions:

\begin{quest}
Suppose that two Heegaard surfaces for a knot-manifold become
isotopic after infinitely many fillings.    Are they equivalent via
annulus swaps?
\end{quest}

If so, then \fullref{lemSwapWGammaPrim} would imply that the
answers to the following is also ``yes'':

\begin{quest}
Suppose that two Heegaard surfaces for a knot-manifold become
isotopic after infinitely many fillings.    Are they weakly
$\gamma$--primitive?
\end{quest}

\bibliographystyle{gtart}
\bibliography{link}

\end{document}